\documentclass[12pt,leqno]{article}
\usepackage{amssymb}
\usepackage[mathscr]{eucal}
\usepackage{amsmath,amssymb,latexsym,theorem,bbm}

\newcommand{\EE}{\mathsf{E}}

\newcommand{\NN}{\mathbb{N}}
\newcommand{\PP}{\mathsf{P}}

\newcommand{\RR}{\mathbb{R}}

\newcommand{\cA}{{\mathcal A}}
\newcommand{\cB}{{\mathcal B}}

\newcommand{\cG}{{\mathcal G}}
\newcommand{\cF}{{\mathcal F}}

\newcommand{\dd}{\mathrm{d}}
\newcommand{\ee}{\mathrm{e}}

\newcommand{\Cov}{\operatorname{Cov}}

\newcommand{\tr}{\operatorname{tr}}

\renewcommand{\mid}{\,|\,}

\renewcommand{\leq}{\leqslant}
\renewcommand{\geq}{\geqslant}

\newcommand{\proofend}{\hfill\mbox{$\Box$}}

\numberwithin{equation}{section}

\theoremstyle{change} \theorembodyfont{\em}
\newtheorem{Lem}{Lemma.}[section]
\newtheorem{Thm}{Theorem.}[section]
\newtheorem{Pro}{Proposition.}[section]

\newtheorem{Def}{Definition.}[section]

\theorembodyfont{\rm}
\newtheorem{Rem}{Remark.}[section]

\begin{document}

\begin{center}
 {\bfseries\Large Representations of multidimensional linear process bridges} \\[5mm]

 {\sc\large M\'aty\'as $\text{Barczy}^{*}$,
            \ Peter $\text{Kern}^{\star,\diamond}$}
\end{center}

\vskip0.2cm

\noindent * Faculty of Informatics, University of Debrecen,
            Pf.~12, H--4010 Debrecen, Hungary.

\noindent $\star$ Mathematical Institute, Heinrich-Heine-Universit\"at D\"usseldorf, 
         Universit\"atsstr.~1, D-40225 D\"usseldorf, Germany.

\noindent e--mails: barczy.matyas@inf.unideb.hu (M. Barczy),
                    kern@math.uni-duesseldorf.de (P. Kern).

\noindent $\diamond$ Corresponding author.

\renewcommand{\thefootnote}{}
\vspace*{0.2cm}
\footnote{M. Barczy has been supported by the Hungarian Scientific Research Fund 
 under Grant No.\ OTKA T-079128.
This work has been finished while M. Barczy was on a post-doctoral position at the Laboratoire 
 de Probabilit\'es et Mod\`{e}les Al\'eatoires, University Pierre-et-Marie Curie, 
 thanks to NKTH-OTKA-EU FP7 (Marie Curie action) co-funded 'MOBILITY' Grant No. MB08-A 81263.}

\vspace*{-10mm}

\begin{abstract}
We derive bridges from general multidimensional linear non time-homogeneous processes using only the transition
 densities of the original process giving their integral representations (in terms of a standard Wiener process) and so-called anticipative representations.
We derive a stochastic differential equation satisfied by the integral representation and we prove a usual
 conditioning property for general multidimensional linear process bridges.
We specialize our results for the one-dimensional case; especially, we study one-dimensional
 Ornstein-Uhlenbeck bridges.
\end{abstract}

\section{Introduction}\label{SECTION_MOTIVATION}

In this paper we deal with deriving bridges from general multidimensional linear processes
 giving their integral representations (in terms of a standard Wiener process) and
 so-called anticipative representations.
Our results are also specialized for the one-dimensional case.
A bridge process is a stochastic process that is pinned to some fixed point at a future
 time point.
Important examples are provided by Wiener bridges, Bessel bridges and general Markovian bridges,
 which have been extensively studied and find numerous applications.
See, for example, Karlin and Taylor \cite[Chapter 15]{KarTay2}, Fitzsimmons, Pitman and Yor
 \cite{FitPitYor}, Privault and Zambrini \cite{PriZam}, Delyon and Hu \cite{DelHu},
 Gasbarra, Sottinen and Valkeila \cite{GasSotVal}, Goldys and Maslowski \cite{GolMas},
  Chaumont and Uribe Bravo \cite{ChaUri} and Baudoin and Nguyen-Ngoc \cite{BauNgu}.
 Recently, Hoyle, Hughston and Macrina \cite{HoyHugMac} studied so-called L\'evy random bridges,
 that is L\'evy processes conditioned to have a prespecified marginal law at the endpoint of the bridge
 (see also the Ph.D. dissertation of Hoyle \cite{Hoy}).
Bichard \cite{Bic} considered so-called bridged Wiener sheets, that is Wiener sheets
 which are forced to take some values along specified curves.

In what follows first we give a motivation for our multidimensional results by presenting different
 representations of the one-dimensional Ornstein-Uhlenbeck bridges, and then we briefly summarize
 the structure of the paper.

\bigskip
\smallskip

\noindent{\bf Motivation: representations of one-dimensional Ornstein-Uhlenbeck bridges}

\smallskip
\smallskip

Let \ $(B_t)_{t\geq 0}$ \ be a standard Wiener process and for
$q\not=0$, $\sigma\not=0$ let us consider the stochastic differential equation (SDE)
\begin{align}\label{OU_egyenlet}
  \begin{cases}
   \dd Z_t=q\,Z_t\,\dd t+\sigma\,\dd B_t,\qquad t\geq0,\\
   \phantom{\dd} Z_0=0.
  \end{cases}
 \end{align}
It is known that there exists a strong solution of this SDE, namely
\begin{align}\label{OU_egyenlet_megoldas}
  Z_t=\sigma\int_0^t\ee^{q(t-s)}\,\dd B_s,\qquad t\geq 0,
 \end{align}
 and strong uniqueness for the SDE \eqref{OU_egyenlet} holds.
 The process \ $(Z_t)_{t\geq 0}$ \ is called a one-dimensional Ornstein-Uhlenbeck
 process (OU-process). It is a time-homogeneous Gauss-Markov process with
 transition densities
 \begin{equation}\label{OU_densities}
  p_t^Z(x,y)
  =\frac{1}{\sqrt{2\pi\sigma^2 \kappa_q(t)}}
     \exp\left\{-\frac{(y-\ee^{qt}x)^2}{2\sigma^2 \kappa_q(t)}\right\},
       \qquad t>0,\quad x,y\in\RR,
 \end{equation}
 where we set
 \begin{equation}\label{kappa_def}
  \kappa_q(t):=\frac{\ee^{2qt}-1}{2q}=\frac{\ee^{qt}}{q}\sinh(qt),\qquad t\geq0.
 \end{equation}
For \ $a,b\in\RR$ \ and \ $T>0,$ \ by an Ornstein-Uhlenbeck bridge from $a$ to $b$
over the time interval \ $[0,T]$ \ derived from \ $Z$ \ we understand a Markov process
\ $(U_t)_{t\in[0,T]}$ \ with initial distribution \ $P(U_0=a)=1,$ \ with
\ $P(U_T=b)=1$ \ and with transition densities
 \begin{equation}\label{bridge_densities}
    p_{s,t}^U(x,y)=\frac{p_{t-s}^Z(x,y)\,p_{T-t}^Z(y,b)}{p_{T-s}^Z(x,b)},
      \quad x,y\in\RR,\quad 0\leq s<t<T.
 \end{equation}
We also note that \ $U_t$ \ converges almost surely to \ $b$ \ as \ $t\uparrow T$,
 \ see, e.g., Fitzsimmons, Pitman and Yor \cite[Proposition 1]{FitPitYor}.
For the construction of bridges derived from a general time-homogeneous
 Markov process using only its transition densities, see, e.g., Barczy and Pap \cite{BarPap}
 and Chaumont and Uribe Bravo \cite{ChaUri}.
Standard calculations yield that for \ $x,y\in\RR$ \ and \ $0\leq s<t<T$,
 \begin{equation}\label{OU_bridge_densities}
   \frac{p_{t-s}^Z(x,y)\,p_{T-t}^Z(y,0)}{p_{T-s}^Z(x,0)}
   =\frac{1}{\sqrt{2\pi\sigma(s,t)}}\,
      \exp\left\{-\frac{\left(y-\frac{\sinh(q(T-t))}{\sinh(q(T-s))}\,x\right)^2}{2\sigma(s,t)}\right\},
 \end{equation}
which is a Gauss density (as a function of \ $y$) with mean \ $\frac{\sinh(q(T-t))}{\sinh(q(T-s))}\,x$ \ and
 variance $\sigma(s,t)$, where for all
 \ $0\leq s\leq t<T,$
 \begin{equation}\label{OU_bridge_var}
   \sigma(s,t):=\sigma^2\,\frac{ \kappa_q(T-t)\kappa_q(t-s)}{\kappa_q(T-s)}
               =\frac{\sigma^2}{q}\,\frac{\sinh(q(T-t))\sinh(q(t-s))}{\sinh(q(T-s))}.
 \end{equation}

Note that if \ $\sigma=0$ \ then for any \ $q\in\RR$ \ the unique
(deterministic) solution of  \eqref{OU_egyenlet}
 is \ $Z_t=0$ \ for all \ $t\geq0$ \ (which coincides with its own bridge from 0 to 0).
 On the other hand, if \ $q=0$ \ and \ $\sigma\not=0$, \ the unique strong solution of
 the SDE \eqref{OU_egyenlet} is the Wiener process \ $Z_t=\sigma B_t,$ $t\geq 0,$ \ and
 it is well known that the Wiener bridge \ $(\widetilde U_t)_{t\in[0,T]}$ \ from $0$ to $0$ over
 \ $[0,T]$ \ derived from \ $Z=\sigma B$ \ admits the (stochastic) integral representation
 \begin{equation}\label{W_bridge}
   \widetilde U_t=\sigma\int_0^t\frac{T-t}{T-s}\,\dd B_s,\qquad t\in[0,T),
 \end{equation}
see, e.g., Section 5.6.B in Karatzas and Shreve \cite{KarShr}.
Moreover, one easily verifies that \ $(\widetilde U_t)_{t\in[0,T]}$
\ is a Markov process with transition densities
 \begin{align}\label{SEGED24}
   p_{s,t}^{\widetilde U}(x,y)
       =\frac{1}{\sqrt{2\pi\widetilde\sigma(s,t)}}\,
        \exp\left\{-\frac{\left(y-\frac{T-t}{T-s}\,x\right)^2}
                 {2\widetilde\sigma(s,t)}\right\},\quad x,y\in\RR,\quad 0\leq s<t<T,
 \end{align}
where \ $\widetilde\sigma(s,t):=\sigma^2\frac{(T-t)(t-s)}{T-s}$ \ for all \ $0\leq s<t<T,$
 \ and that \eqref{bridge_densities} is satisfied with \ $b=0$, \ $U$ \ being replaced by
  \ $\widetilde U$ \  and
 \[
    p_t^Z(x,y)
       =\frac{1}{\sqrt{2\pi t\sigma^2}}\,
        \exp\left\{-\frac{(y-x)^2}{2t\sigma^2}\right\},\quad x,y\in\RR,\quad t>0.
 \]
Comparing \eqref{OU_bridge_densities} with \eqref{SEGED24}, it is quite
reasonable that an integral representation for the Ornstein-Uhlenbeck bridge
from $0$ to $0$ over \ $[0,T]$ \ derived from the process \ $Z$ \ given by the
SDE \eqref{OU_egyenlet} should have the form
 $$
   U_t=\sigma\int_0^t\frac{\sinh(q(T-t))}{\sinh(q(T-s))}\,\dd B_s,
      \qquad t\in[0,T),
 $$
and in fact this is made precise in the sequel. We will further
consider general multivariate linear process bridges.

Besides the integral representation \eqref{W_bridge} of the Wiener
bridge \ $(\widetilde U_t)_{t\in[0,T]}$ \ from \ $0$ \ to \ $0$ \ over
 \ $[0,T],$ \ one can find two equivalent representations in the literature.
Namely, by Section 5.6.B in Karatzas and Shreve \cite{KarShr},
  \begin{align}\label{W_bridge2}
    \begin{cases}
      \dd \widetilde U_t=-\frac{1}{T-t}\widetilde U_t\,\dd t+\dd B_t,\qquad t\in[0,T),\\
      \phantom{\dd} \widetilde U_0=0,
    \end{cases}
 \end{align}
and
\begin{align}\label{W_bridge3}
  \widehat{U_t}=B_t-\frac{t}{T}B_T,\qquad t\in[0,T].
\end{align}
The representation \eqref{W_bridge} with \ $\sigma=1$ \ is just a strong solution of
 the SDE \eqref{W_bridge2}.
 So, the equations \eqref{W_bridge} with \ $\sigma=1$ \ and \eqref{W_bridge2} define the same process
 \ $(\widetilde U_t)_{t\in[0,T]}.$ \ However, the equation \eqref{W_bridge3} does not define
  the same process as the equations \eqref{W_bridge} with \ $\sigma=1$ \ and \eqref{W_bridge2}. The equality
  between representations \eqref{W_bridge} with \ $\sigma=1,$ \eqref{W_bridge2} and \eqref{W_bridge3} is only
  an equality in law, i.e., they determine the same probability measure on \ $(C([0,T]),\cB(C([0,T])))$, \ where
  \ $C([0,T])$ \ denotes the set of all real-valued continuous functions on \ $[0,T]$ \ and
  \ $\cB(C([0,T]))$ \ is the Borel $\sigma$-algebra on it.
The fact that the processes \ $\widetilde U$ \ and \ $\widehat U$ \ are different
  follows from the fact that the process \ $\widetilde U$ \ is adapted to the filtration generated
  by \ $B,$ \ while the process \ $\widehat U$ \ is not. Indeed, to construct \ $\widehat U$ \ we need
  the random variable \ $B_T.$ \ One can call \eqref{W_bridge3} a non-adapted, anticipative representation
  of a Wiener bridge. The attribute anticipative indicates that
  for the definition of \ $\widehat{U_t}$ \ we use the random
  variable \ $B_T,$ \ where the time point \ $T$ \ is after the time point \
  $t.$ \

A similar anticipative representation of an
  Ornstein-Uhlenbeck bridge derived from the SDE \eqref{OU_egyenlet}
  can be found on page 378 in Donati-Martin \cite{DonMar} and
  in Lemma 1 in Papie\.{z} and Sandison \cite{PapSan}.
  Donati-Martin gave an anticipative representation of an Ornstein-Uhlenbeck
  bridge from $a=0$ to $b=0$ derived from the SDE \eqref{OU_egyenlet} with $q<0$ and $\sigma=1$, while
  Papie\.{z} and Sandison formulated their lemma in case of arbitrary starting point $a$
  and ending point $b,$ \ but only
  for special values of \ $q$ \ and \ $\sigma,$ \ but their proof is  also valid
  for all \ $q\ne 0$ \ and \ $\sigma\ne 0$ \ (see our Remark \ref{Remark7}).

 Moreover, concerning the relationship between Wiener
processes and Wiener bridges, by Problem 5.6.13 in Karatzas and
Shreve \cite{KarShr}, if \ $T>0$ \ is fixed and \ $(B_t)_{t\geq 0}$
\  is a standard Wiener process (starting from $0$), then for all
  \ $n\in\NN,$ \ $0<t_1<\ldots<t_n<T,$ \ the conditional distribution of
 \ $(B_{t_1},\ldots,B_{t_n})$ \ given \ $B_T=0$ \ equals the distribution
of \ $(\widetilde U_{t_1},\ldots,\widetilde U_{t_n}),$ \ where \ $\widetilde U$ \
is given by \eqref{W_bridge} with \ $\sigma=1$ \ or by \eqref{W_bridge2}.

Finally, we note that the transition densities \ $p_{s,t}^U(x,y)$, $x,y\in\RR$, $0\leq s<t<T$,
 \ of the process bridge \ $(U_t)_{t\in[0,T]}$ \ can be derived using Doob's \ $h$-transform
 (see Doob \cite{Doob}).
In Section \ref{SECTION_MULTI} we briefly study this approach for general multivariate linear process bridges.

\bigskip
\smallskip

\noindent{\bf Structure of the paper}

\smallskip
\smallskip

In Section \ref{SECTION_MULTI} we derive multidimensional linear process bridges from a multidimensional
 linear non time-homogeneous process \ $Z$ \ given by the SDE \eqref{gen_multivariate_system} using only the
 transition densities of \ $Z$, \ see Theorem \ref{gen_multivariate_intrep_bridge} and
 Definition \ref{DEF_bridge_multidim}.
We also give an integral and a so-called anticipative representation of the derived bridge,
 see formulae \eqref{gen_multivariate_intrep} and \eqref{anticipative_def}, respectively.
We derive an SDE satisfied by this integral representation, see Theorem \ref{LEMMA_multi_DE},
 and in Proposition \ref{bridge_conditioning} we prove a usual conditioning property for general
 multidimensional linear process bridges.
 In Remark \ref{Remark18} we point out that the integral representation and anticipative representation
 of the bridge are quite different.
To shed more light on the different behavior of the different bridge representations,
 in a companion paper we will study sample path deviations of the Wiener
 process and the Ornstein-Uhlenbeck process from its bridges, see Barczy and Kern \cite{BarKer}.
In Remark \ref{Remark14} we study that the SDE derived for the integral representation can be considered
 as a consequence of Proposition 3 in Delyon and Hu \cite{DelHu}.
We use the expression 'can be considered' since the definition of bridges given in Delyon and Hu \cite{DelHu}
 and in the present paper are different.
We have a different approach coming from the possibility that in our special case we are able
 to explicitly calculate the transition densities of the bridge from which we deduce an integral
 representation and finally end up with the same SDE of Proposition 3 in Delyon and Hu \cite{DelHu}
 such that this integral representation is a strong solution of the above mentioned SDE.
We also note that the SDE of Proposition 3 in Delyon and Hu \cite{DelHu} contains the solution of a
 deterministic differential equation (see Remark \ref{Remark14}) which solution always remains abstract,
 while in our special case we have an explicit solution via evolution matrices
 (see Section \ref{SECTION_MULTI}).
Theorem 2 of Delyon and Hu \cite{DelHu} gives anticipative representations of general
 multidimensional conditioned diffusions, in Remark \ref{Remark12} we compare our results for anticipative
 representations of process bridges in the multidimensional case with the corresponding results of
 Delyon and Hu \cite{DelHu}.
Concerning anticipative representations of one-dimensional Gauss bridges see the recent paper
 of Gasbarra, Sottinen and Valkeila \cite{GasSotVal}.
Remark \ref{Remark17} is devoted to studying the construction of a bridge by subtracting from
 a process its conditional expectation given the process at a prescribed time point (the endpoint of the
 bridge).

In Section \ref{SECTION_ONEDIMENSIONAL} we formulate our multidimensional results in case of dimension one
 which includes also the study of usual Ornstein-Uhlenbeck bridges.
We note that not all of the results are immediate consequences of the multidimensional ones
 and in case of dimension one we can give an illuminating explanation for the anticipative representation
 motivated by Lemma 1 in Papie\.{z} and Sandison \cite{PapSan}, see Remark \ref{Remark7}.
Moreover, in Remark \ref{Remark15} we discuss the connections between Propositions 4 and 9 in Gasbarra,
 Sottinen and Valkeila \cite{GasSotVal} and our Theorems \ref{PROPOSITION6} and \ref{THEOREM2}
 (anticipative and integral representation of the bridge in case of dimension one).

Section \ref{SECTION_PROOFS} (Appendix) is devoted to the proofs of our results and
 contains also a supplement for our assumption on Kalman type matrices
 (introduced in Section \ref{SECTION_MULTI}).

\section{Multidimensional linear process bridges}\label{SECTION_MULTI}

Let \ $\NN$, \ $\RR$ \ and \ $\RR_+$ \ denote the set of positive integers, real numbers and
 non-negative real numbers, respectively.
For all \ $n,m\in\NN$, \ let \ $\RR^{n\times m}$ \ and \ $I_n$ \ denote the set of \ $n\times m$
 \ matrices with real entries and the \ $n\times n$ \ identity matrix, respectively.

For all \ $d,p\in\NN$, \ let us consider a general $d$-dimensional linear process given
 by the linear SDE
\begin{align}\label{gen_multivariate_system}
{\rm d}\mathbf Z_{t}
  =\big(Q(t)\mathbf Z_{t}+\mathbf r(t)\big)\,{\rm d} t
    +\Sigma(t)\,{\rm d}\mathbf B_t,\qquad t\geq0,
\end{align}
with continuous functions $Q:\RR_+\to\RR^{d\times d}$, $\Sigma:\RR_+\to\RR^{d\times p}$ and
 $\mathbf r:\RR_+\to\RR^d$, where $(\mathbf B_t)_{t\geq 0}$ is a $p$-dimensional standard Wiener process
 on a filtered probability space \ $(\Omega,\cF,(\cF_t)_{t\geq 0},P)$ \ satisfying the usual conditions
 (the filtration being constructed by the help of \ $\mathbf B$),
  i.e., \ $(\Omega,\cF,P)$ \ is complete,
 \ $(\cF_t)_{t\geq 0}$ \ is right continuous, \ $\cF_0$ \ contains all the $P$-null sets in
 \ $\cF$ \ and \ $\cF_\infty=\cF$, \ where \ $\cF_\infty:=\sigma\left(\bigcup_{t\geq 0}\cF_t\right)$,
 see, e.g., Karatzas and Shreve \cite[Section 5.2.A]{KarShr}.
It is known that there exists a strong solution of the SDE
\eqref{gen_multivariate_system}, namely
\begin{align}\label{gen_multivariate_int_rep}
\mathbf Z_{t}=\Phi(t)\left[\mathbf Z_{0}+\int_0^t\Phi^{-1}(s)\mathbf r(s)\,{\rm d} s+\int_0^t\Phi^{-1}(s)\Sigma(s)\,{\rm d}\mathbf B_s\right],\qquad t\geq 0,
\end{align}
where $\Phi$ is a solution to the deterministic matrix differential equation
 $\Phi'(t)=Q(t)\Phi(t)$, $t\geq 0$, with $\Phi(0)=I_d$ \ and strong uniqueness
 for the SDE \eqref{gen_multivariate_system} holds, see,
 e.g., Karatzas and Shreve \cite[Section 5.6]{KarShr}.
The unique solution of the above matrix differential equation can be given as $\Phi(t)=E(t,0)$, $t\geq 0$,
 in terms of the evolution matrices
\begin{align*}
E(t,s)=I_d + \int_{s}^tQ(t_{1})\,{\rm d} t_{1}
       +\sum_{k=2}^\infty\int_{s}^t\int_{s}^{t_1}\cdots\int_{s}^{t_{k-1}}Q(t_{1})
          \cdots Q(t_{k})\,{\rm d} t_{k}{\rm d} t_{k-1}\cdots{\rm d} t_{1}
\end{align*}
for $s,t\geq0$. Indeed, by Theorem 1.8.2 in Conti \cite{Con}, the general $d$-dimensional solution
 $\mathbf y(t):\RR_{+}\to\RR^{d}$ of $\mathbf y'(t)=Q(t)\mathbf y(t),$ $t\geq 0$, is represented by
 $\mathbf y(t)=E(t,s)\mathbf y(s)$ for all $s,t\geq0$, which shows that $\Phi(t)=E(t,0)$, $t\geq 0$.
Note that, since $Q$ is continuous, there exists an $L>0$ such that $\Vert Q(u)\Vert\leq L$
 for all $u\in[\min(s,t),\max(s,t)]$, $s,t\geq 0$
 (with some fixed matrix norm $\Vert .\Vert$ on $\RR^{d\times d}$), and hence one easily calculates
 $\Vert E(t,s)\Vert\leq{\rm e}^{L|t-s|}$. Note also that if $Q(t)= Q\in\RR^{d\times d}$, $t\geq 0$,
 is constant then $E(t,s)={\rm e}^{(t-s)Q}$ for $t,s\geq 0$, and
 hence $\Phi(t)={\rm e}^{t Q},$ $t\geq 0.$

We will make frequent use of the following properties of evolution matrices stated as equations
 (1.9.2) and (1.9.3) in Conti \cite{Con}. For all $r,s,t\geq0$ we have
\begin{align}
& E(t,s)E(s,r)=E(t,r),\label{evolution}\\
& E(t,t)=I_d,\quad\quad E(t,s)^{-1}=E(s,t),\label{evoinvert}\\
& \partial_{1}E(t,s)=Q(t)E(t,s),\quad\quad\partial_{2}E(t,s)=-E(t,s)Q(s).
\label{evoderive}
\end{align}

The unique strong solution of the SDE
\eqref{gen_multivariate_system} can now be written as
$$
 \mathbf Z_{t}
  =E(t,0)\mathbf Z_{0}+\int_0^tE(t,s)\mathbf r(s)\,{\rm d} s+\int_0^tE(t,s)\Sigma(s)\,{\rm d}\mathbf B_s,
    \qquad t\geq 0.
$$
Here and in what follows we assume that $\mathbf Z_{0}$ has a Gauss distribution independent
 of the Wiener process $(\mathbf B_{t})_{t\geq0}$.
Then we may define the filtration \ $(\cF_t)_{t\geq 0}$ \ such that
 \ $\sigma\{\mathbf Z_0,\mathbf B_s:\,0\leq s\leq t\}\subset\cF_t$ for all $t\geq 0$,
 see, e.g., Karatzas and Shreve \cite[Section 5.2.A]{KarShr}.

We will call the process $(\mathbf Z_{t})_{t\geq0}$ a $d$-dimensional linear process.

One can easily derive that for $0\leq s\leq t$ we have
\begin{equation}\label{increments_solution}
 \mathbf Z_{t}=E(t,s)\mathbf Z_{s}+\int_s^tE(t,u)\mathbf r(u)\,{\rm d} u
   + \int_{s}^{t}E(t,u)\Sigma(u)\,{\rm d}\mathbf B_{u}.
\end{equation}
Hence, given $\mathbf Z_s=\mathbf x,$ the distribution of $\mathbf Z_t$ does not depend on
 $(\mathbf Z_{u})_{u\in[0,s)}$ and thus $(\mathbf Z_{t})_{t\geq0}$ is a Gauss-Markov process
 (see, e.g., Karatzas and Shreve \cite[Problem 5.6.2]{KarShr}).
For any $0\leq s\leq t$ and $\mathbf x\in\RR^d$ let us define
$$
 \mathbf m_{\mathbf x}^+(s,t)
   :=\mathbf x+\int_{s}^tE(s,u)\mathbf r(u)\,{\rm d} u
     \quad\text{ and }\quad
 \mathbf m_{\mathbf x}^-(s,t):=\mathbf x-\int_{s}^tE(t,u)\mathbf r(u)\,{\rm d}u.
$$
Then for any $\mathbf x\in\RR^d$ and $0\leq s<t$ the conditional distribution of $\mathbf Z_t$ given
 $\mathbf Z_s=\mathbf x$ is Gauss with mean
 $$
   \mathbf m_{\mathbf x}(s,t):=E(t,s)\mathbf m_{\mathbf x}^+(s,t)
       =E(t,s)\mathbf x+\int_s^tE(t,u)\mathbf r(u){\rm d}u,
 $$
and with covariance matrix of Kalman type (see Kalman \cite{Kal})
\begin{equation*}
\kappa(s,t):= \int_{s}^tE(t,u)\Sigma(u)\Sigma(u)^\top E(t,u)^\top\,{\rm d}u.
\end{equation*}
The matrices $\kappa(s,t)$ are symmetric and positive semi-definite for all $0\leq s<t$,
 and in what follows we put the following assumption:
 \begin{equation}\label{kappa_assumption}
  \kappa(s,t)\text{ is positive definite for all }0\leq s<t.
 \end{equation}
From control theory of linear systems we owe sufficient conditions for positive
 definiteness of the Kalman matrices (see, e.g., Theorems 7.7.1 - 7.7.3 in Conti \cite{Con})
 which we present in the Appendix, see Proposition \ref{PROPOSITION_kappa}.

Hence the transition densities of the Gauss-Markov process $(\mathbf Z_{t})_{t\geq0}$ read as
\begin{align}\label{gen_multivariate_density}
p_{s,t}^{\mathbf Z}(\mathbf x,\mathbf y)
 =\frac{1}{\sqrt{(2\pi)^d\det\kappa(s,t)}}
    \exp\left\{-\frac{1}{2}\big\langle\kappa(s,t)^{-1}(\mathbf y-\mathbf m_{\mathbf x}(s,t)),
             \mathbf y-\mathbf m_{\mathbf x}(s,t)\big\rangle\right\}
\end{align}
for all $0\leq s<t$ and $\mathbf x,\mathbf y\in\mathbb R^d$. Our aim is to derive
 a process bridge from $\mathbf Z$, namely, we will consider a bridge from $\mathbf a$
 to $\mathbf b$ over the time interval $[0,T]$, where $\mathbf a,\mathbf b\in\RR^d$ and $T>0$.
Generalizing the formula (2.7) in Fitzsimmons, Pitman and Yor \cite{FitPitYor} to multidimensional
 non time-homogeneous Markov processes, for fixed $T>0$ we are looking for a Markov process
 $(\mathbf U_t)_{t\in[0,T]}$ with initial distribution $P(\mathbf U_0=\mathbf a)=1$ and
 with transition densities
\begin{equation}\label{gen_multivariate_bridge_densities}
 p_{s,t}^{\mathbf U}(\mathbf x,\mathbf y)
   =\frac{p_{s,t}^{\mathbf Z}(\mathbf x,\mathbf y)
    \,p_{t,T}^{\mathbf Z}(\mathbf y,\mathbf b)}{p_{s,T}^{\mathbf Z}(\mathbf x,\mathbf b)},
    \quad\mathbf x,\mathbf y\in\RR^d,\quad 0\leq s<t<T,
\end{equation}
provided that such a process exists.
To properly speak of $(\mathbf U_t)_{t\in[0,T]}$ as a process bridge, we shall study
 the limit behavior of $\mathbf U_t$ as $t\uparrow T$, namely, we shall show that
 $\mathbf U_{t}\to\mathbf b=:\mathbf U_{T}$ almost surely and also in $L^2$ as $t\uparrow T$
 (see Theorem \ref{gen_multivariate_intrep_bridge}).

Our approach can also be seen in the context of Doob's $h$-transform (see Doob \cite{Doob}) as follows.
For bounded Borel-measurable functions \ $f:\RR_+\times\RR^d\to\RR$ \ one can define a family
 of operators \ $(P_{s,t})_{0\leq s<t}$ \ by
$$
 P_{s,t}f(s,\mathbf x):=\int_{\RR^d}f(t,\mathbf y)
                              \,p_{s,t}^{\mathbf Z}(\mathbf x,\mathbf y)\,{\rm d}\mathbf y
 $$
for $0\leq s<t$ and $\mathbf x\in\RR^d$.
 Then
 \begin{align*}
       & \vert P_{s,t}f(s,\mathbf x) \vert
            \leq \int_{\RR^d}\vert f(t,\mathbf y)\vert
               \,p_{s,t}^{\mathbf Z}(\mathbf x,\mathbf y)\,{\rm d}\mathbf y
            \leq \sup_{\mathbf y\in\RR^d}\vert f(t,\mathbf y)\vert
            <\infty,\\
       & P_{s,t}f(s,\mathbf Z_s) = \EE(f(t,\mathbf Z_t)\mid\mathbf Z_s) \qquad \text{$P$-a.s.},
 \end{align*}
 and the family \ $(P_{s,t})_{0\leq s<t}$ \ forms a hemigroup of transition operators
 for the Markov process \ $\mathbf Z$. \
Indeed, for \ $0\leq s<r<t$ \ and \ $\mathbf x\in\RR^d$ \ we observe
\begin{align*}
 P_{s,r}P_{r,t}f(s,\mathbf x)
  & =\int_{\RR^d}P_{r,t}f(r,\mathbf y)\,p_{s,r}^{\mathbf Z}(\mathbf x,\mathbf y)\,{\rm d}\mathbf y
    =\int_{\RR^d}\int_{\RR^d}f(t,\mathbf z)\,p_{r,t}^{\mathbf Z}(\mathbf y,\mathbf z)\,{\rm d}
     \mathbf z\,p_{s,r}^{\mathbf Z}(\mathbf x,\mathbf y)\,{\rm d}\mathbf y\\
  & =\int_{\RR^d}f(t,\mathbf z)\int_{\RR^d}p_{s,r}^{\mathbf Z}(\mathbf x,\mathbf y)\,p_{r,t}^{\mathbf Z}
      (\mathbf y,\mathbf z)\,{\rm d}\mathbf y\,{\rm d}\mathbf z\\
 & =\int_{\RR^d}f(t,\mathbf z)\,p_{s,t}^{\mathbf Z}(\mathbf x,\mathbf z)\,{\rm d}\mathbf z
   =P_{s,t}f(s,\mathbf x).
\end{align*}
For fixed \ $T>0$ \ and  \ $\mathbf b\in\RR^d$ \ we now define the function
$$
 h:[0,T)\times\RR^d\to\RR_+\quad\text{ by }\quad h(t,\mathbf x)=p_{t,T}^{\mathbf Z}(\mathbf x,\mathbf b),
  \qquad t\in[0,T),\;\; \mathbf x\in\RR^d.
$$
By \eqref{gen_multivariate_density}, \ $h$ \ is positive and bounded on
 \ $[0,t]\times\RR^d$ \ for every \ $0<t<T$.
\ Indeed, \eqref{kappa_assumption} yields that
 \[
    \inf_{s\in[0,t]}\det \kappa(s,T) >0,\qquad t\in[0,T),
 \]
 and hence
 \[
   \sup_{(s,x)\in[0,t]\times\RR^d} \vert h(s,x)\vert
      \leq \left( (2\pi)^d \inf_{s\in[0,t]}\det \kappa(s,T) \right)^{-1/2}<\infty,
      \qquad t\in[0,T).
 \]
This yields that $P_{s,t}h(s,\mathbf x)$ is defined for all $0\leq s<t<T$ and
 $\mathbf x\in\RR^d$, although it can happen that $h$ is not bounded on $[0,T)\times\RR^d$
 (as it is in the case of $\mathbf Z$ being a one-dimensional standard Wiener process).
Then \ $h$ \ is space-time harmonic for the Markov process \ $\mathbf Z$ \ in the sense that
 $$
   P_{s,t}h(s,\mathbf x)
     =\int_{\RR^d}h(t,\mathbf y)\,p_{s,t}^{\mathbf Z}(\mathbf x,\mathbf y)\,{\rm d}\mathbf y
     =\int_{\RR^d}p_{s,t}^{\mathbf Z}(\mathbf x,\mathbf y)\,p_{t,T}^{\mathbf Z}(\mathbf y,\mathbf b)
      \,{\rm d}\mathbf y=p_{s,T}^{\mathbf Z}(\mathbf x,\mathbf b)=h(s,\mathbf x)
 $$
for $0\leq s<t<T$ and $\mathbf x\in\RR^d$.
Now a generalization of Doob's $h$-transform approach (see Doob \cite{Doob}
 gives a new operator hemigroup
 $$
  \widetilde P_{s,t}f=\frac1h\,P_{s,t}(hf),\quad0\leq s<t<T
 $$
with
\begin{align*}
 \widetilde P_{s,t}f(s,\mathbf x) & =\frac1{h(s,\mathbf x)}\,P_{s,t}(hf)(s,\mathbf x)
  =\frac1{h(s,\mathbf x)}\int_{\RR^d}h(t,\mathbf y)f(t,\mathbf y)\,p_{s,t}^{\mathbf Z}
   (\mathbf x,\mathbf y)\,{\rm d}\mathbf y\\
& =\int_{\RR^d}f(t,\mathbf y)\,\frac{p_{s,t}^{\mathbf Z}(\mathbf x,\mathbf y)\,p_{t,T}^{\mathbf Z}
   (\mathbf y,\mathbf b)}{p_{s,T}^{\mathbf Z}(\mathbf x,\mathbf b)}\,{\rm d}\mathbf y
    =\int_{\RR^d}f(t,\mathbf y)\,p_{s,t}^{\mathbf U}(\mathbf x,\mathbf y)\,{\rm d}\mathbf y,
\end{align*}
 where \ $f:\RR_+\times \RR^d\to\RR$ \ is a bounded Borel-measurable function and
  $\mathbf x\in\RR^d$, i.e., the transition operators
  \ $(\widetilde P_{s,t})_{0\leq s<t<T}$ \ belong to a new Markov process \ $(\mathbf U_t)_{0\leq t<T}$,
  \ the desired process bridge, with transition densities \ $(p_{s,t}^{\mathbf U})_{0\leq s<t<T}$ \ given by \eqref{gen_multivariate_bridge_densities}.

For $T>0$, $0\leq s<t<T$ and $\mathbf a,\mathbf b\in\RR^d$, let us define
 \begin{align*}
  &\Gamma(s,t):=E(s,t)\kappa(s,t) = \int_s^t E(s,u)\Sigma(u)\Sigma(u)^\top E(t,u)^\top\,\dd u,\\
  &\Sigma(s,t):=\Gamma(t,T)\Gamma(s,T)^{-1}\Gamma(s,t),
 \end{align*}
 and
 \begin{equation}\label{bridge_mean}
 \mathbf n_{\mathbf a,\mathbf b}(s,t)
  :=\Gamma(t,T)\Gamma(s,T)^{-1}\mathbf m_{\mathbf a}^+(s,t)
     +\Gamma(s,t)^\top\big(\Gamma(s,T)^\top\big)^{-1}\mathbf m_{\mathbf b}^-(t,T).
\end{equation}

In what follows we prove the existence of a Markov process $(\mathbf U_t)_{t\in[0,T]}$ with initial
 distribution $P(\mathbf U_0=\mathbf a)=1$ and with transition densities \ $p_{s,t}^{\mathbf U}$ \ given
 in \eqref{gen_multivariate_bridge_densities} such that $\mathbf U_{t}\to\mathbf b=:\mathbf U_{T}$ almost surely
 and also in $L^2$ as $t\uparrow T$.
First we present an auxiliary lemma.

\begin{Lem}\label{gen_multivariate_bridge}
Let us suppose that condition \eqref{kappa_assumption} holds.
Let $\mathbf b\in\RR^d$ and $T>0$ be fixed. Then for all $0\leq s<t<T$ and $\mathbf x,\mathbf y\in\RR^d$ we have
\begin{align*}
& \frac{p_{s,t}^{\mathbf Z}(\mathbf x,\mathbf y)\,
     p_{t,T}^{\mathbf Z}(\mathbf y,\mathbf b)}{p_{s,T}^{\mathbf Z}(\mathbf x,\mathbf b)}\\
& \quad=\frac1{\sqrt{(2\pi)^d\det\Sigma(s,t)}}\,
  \exp\left\{-\frac12\,\Big\langle\Sigma(s,t)^{-1}\big(\mathbf y-\mathbf n_{\mathbf x,\mathbf b}(s,t)\big),
            \mathbf y-\mathbf n_{\mathbf x,\mathbf b}(s,t)\Big\rangle\right\},
\end{align*}
which is a Gauss density with mean vector $\mathbf n_{\mathbf x,\mathbf b}(s,t)$ and
 with covariance matrix $\Sigma(s,t)$.
\end{Lem}

The proof of Lemma \ref{gen_multivariate_bridge} can be found in the Appendix.

\begin{Thm}\label{gen_multivariate_intrep_bridge}
Let us suppose that condition \eqref{kappa_assumption} holds.
For fixed $\mathbf a,\mathbf b\in\RR^d$ and $T>0$, let the process
 $(\mathbf U_{t})_{t\in[0,T)}$ be given by
 \begin{equation}\label{gen_multivariate_intrep}
  \mathbf U_{t}
  := \mathbf n_{\mathbf a,\mathbf b}(0,t)
     +\Gamma(t,T)\int_{0}^t\Gamma(u,T)^{-1}\Sigma(u)\,{\rm d}\mathbf B_{u},
     \quad t\in[0,T).
 \end{equation}
Then for any $t\in[0,T)$ the distribution of $\mathbf U_t$ is Gauss with mean
 $\mathbf n_{\mathbf a,\mathbf b}(0,t)$ and covariance matrix $\Sigma(0,t)$.
Especially, $\mathbf U_{t}\to\mathbf b$ almost surely (and hence in probability) and
 in $L^2$ as $t\uparrow T$.
Hence the process $(\mathbf U_{t})_{t\in[0,T)}$ can be extended to an almost surely (and hence stochastically)
 and \ $L^2$-continuous process $(\mathbf U_{t})_{t\in[0,T]}$ with $\mathbf U_{0}=\mathbf a$ and
 $\mathbf U_{T}=\mathbf b$.
Moreover, $(\mathbf U_{t})_{t\in[0,T]}$ is a Gauss-Markov process and for any
 $\mathbf x\in\RR^d$ and $0\leq s<t<T$ the transition density
 $\RR^d\ni \mathbf y\mapsto p_{s,t}^{\mathbf U}(\mathbf x,\mathbf y)$ of $\mathbf U_{t}$ given
 $\mathbf U_{s}=\mathbf x$ is given by
 \[
    p_{s,t}^{\mathbf U}(\mathbf x,\mathbf y)
      = \frac1{\sqrt{(2\pi)^d\det\Sigma(s,t)}}\,
  \exp\left\{-\frac12\,\Big\langle\Sigma(s,t)^{-1}\big(\mathbf y-\mathbf n_{\mathbf x,\mathbf b}(s,t)\big),
            \mathbf y-\mathbf n_{\mathbf x,\mathbf b}(s,t)\Big\rangle\right\},
 \]
 which coincides with the density given in Lemma \ref{gen_multivariate_bridge}.
\end{Thm}

The proof of Theorem \ref{gen_multivariate_intrep_bridge} can be found in the Appendix.

\begin{Def}\label{DEF_bridge_multidim}
Let $(\mathbf Z_{t})_{t\geq0}$ be the $d$-dimensional linear process given by the SDE
 \eqref{gen_multivariate_system} with an initial Gauss random variable $\mathbf Z_0$ independent
 of $(\mathbf B_{t})_{t\geq0}$
and let us assume that condition \eqref{kappa_assumption} holds.
For fixed $\mathbf a,\mathbf b\in\RR^d$ and $T>0$, the process $(\mathbf U_{t})_{t\in[0,T]}$
 defined in Theorem \ref{gen_multivariate_intrep_bridge} is called a linear process bridge
 from $\mathbf a$ to $\mathbf b$ over $[0,T]$ derived from \ $\mathbf Z$.
More generally, we call any almost surely continuous (Gauss) process on the time interval
 $[0,T]$ having the same finite-dimensional distributions as $(\mathbf U_{t})_{t\in[0,T]}$
 a multidimensional linear process bridge from $\mathbf a$ to $\mathbf b$ over $[0,T]$ derived
 from \ $\mathbf Z$.
\end{Def}

Formula \eqref{gen_multivariate_intrep} can be considered as an integral representation of
 the linear process bridge \ $\mathbf U$.

In the next theorem we present an SDE satisfied by the linear process bridge \ $\mathbf U$.

\begin{Thm}\label{LEMMA_multi_DE}
Let us suppose that condition \eqref{kappa_assumption} holds.
The process $(\mathbf U_{t})_{t\in[0,T)}$ defined by
\eqref{gen_multivariate_intrep} is a strong solution of the
linear SDE
\begin{equation}\label{gen_mult_bridge_sde}
 \begin{split}
{\rm d}\mathbf U_{t}
 =&\Big[\big(Q(t)-\Sigma(t)\Sigma(t)^\top E(T,t)^\top\Gamma(t,T)^{-1}\big)\mathbf U_{t}\\
  &+\Sigma(t)\Sigma(t)^\top\big(\Gamma(t,T)^\top\big)^{-1}
  \mathbf m_{\mathbf b}^-(t,T)
    +\mathbf r(t)\Big]{\rm d} t+\Sigma(t)\,{\rm d}\mathbf B_{t}
\end{split}
\end{equation}
for $t\in[0,T)$ and with initial condition $\mathbf U_{0}=\mathbf a$, and strong uniqueness
 for the SDE \eqref{gen_mult_bridge_sde} holds.
\end{Thm}

The proof of Theorem \ref{LEMMA_multi_DE} can be found in the Appendix.

Now we turn to give alternative representations of the bridge.
The next theorem is about the existence of a so-called anticipative representation of the
 bridge which is a weak solution to the bridge SDE \eqref{gen_mult_bridge_sde}.

\begin{Thm}\label{anticipative_bridge}
Let $\mathbf a,\mathbf b\in\RR^d$ and $T>0$ be fixed. Let $(\mathbf Z_t)_{t\geq 0}$
 be the linear process given by the SDE \eqref{gen_multivariate_system} with initial condition
 $\mathbf Z_0=\mathbf 0$ and let us assume that condition \eqref{kappa_assumption} holds.
Then the process $(\mathbf Y_t)_{t\in[0,T]}$ given by
\begin{equation}\label{anticipative_def}
 \mathbf Y_t:=\Gamma(t,T)\Gamma(0,T)^{-1}\mathbf a
  + \mathbf Z_t-\Gamma(0,t)^\top\big(\Gamma(0,T)^\top\big)^{-1}(\mathbf Z_T-\mathbf b),\quad t\in[0,T],
\end{equation}
 equals in law the linear process bridge $(\mathbf U_t)_{t\in[0,T]}$
 from $\mathbf a$ to $\mathbf b$ over $[0,T]$
 derived from $\mathbf Z$.
\end{Thm}

The proof of Theorem \ref{anticipative_bridge} can be found in the Appendix.

Next we present a usual conditioning property for multidimensional linear process bridges.

\begin{Pro}\label{bridge_conditioning}
Let $\mathbf a,\mathbf b\in\RR^d$ and $T>0$ be fixed. Let $(\mathbf Z_t)_{t\geq 0}$ be the $d$-dimensional
 linear process given by the SDE \eqref{gen_multivariate_system} with initial condition $\mathbf Z_0=\mathbf a$
 and let us assume that condition \eqref{kappa_assumption} holds.
 Let $n\in\NN$ and $0<t_1<t_2<\ldots<t_n<T$. Then the conditional distribution of
 $(\mathbf Z_{t_1}^\top,\ldots,\mathbf Z_{t_n}^\top)^\top$ given $\mathbf Z_T=\mathbf b$ equals
 the distribution of $(\mathbf U_{t_1}^\top,\ldots,\mathbf U_{t_n}^\top)^\top,$ where
 $(\mathbf U_t)_{t\in[0,T]}$ is the linear process bridge from $\mathbf a$ to $\mathbf b$ over $[0,T]$
 derived from $(\mathbf Z_t)_{t\geq 0}$.
\end{Pro}

The proof of Proposition \ref{bridge_conditioning} can be found in the Appendix.
One can also realize that in case of time-homogeneity Proposition \ref{bridge_conditioning}
 is a consequence of Proposition 1 in Fitzsimmons, Pitman and Yor \cite{FitPitYor}.
To be more precise, restricting considerations to our situation of a $d$-dimensional linear process
 $(\mathbf Z_t)_{t\geq 0}$ given by the SDE \eqref{gen_multivariate_system} with initial condition
 $\mathbf Z_0=\mathbf a$, Proposition 1 in Fitzsimmons, Pitman and Yor \cite{FitPitYor} states that
 if $\mathbf Z$ is time-homogeneous (with transition densities
 $p_{t}^{\mathbf Z}(\mathbf x,\mathbf y):=p_{s,s+t}^{\mathbf Z}(\mathbf x,\mathbf y)$ for all $s,t\geq0$
 and $\mathbf x,\mathbf y\in\RR^d$), then for fixed $\mathbf a,\mathbf b\in\RR^d$, $T>0$ there exists a unique
 probability measure $\widetilde\PP_{\mathbf a,\mathbf b}^T$ on
 $(\widetilde\Omega, \widetilde {\mathcal F}_{T-})$
 such that $({\mathbf Z}_t)_{t\in[0,T)}$ under $\widetilde\PP_{\mathbf a,\mathbf b}^T$
 is a (non time-homogeneous) Markov process with transition densities given
 by \eqref{gen_multivariate_bridge_densities},
 where $\widetilde\Omega$ is the set of all real-valued c\`adl\`ag functions on $[0,\infty)$,
 $(\widetilde\cF_t)_{t\geq 0}$ is the natural (uncompleted) filtration of the coordinate process
 $(\mathbf Z_t)_{t\geq 0}$ on $\widetilde\Omega$ (which we also denote by $\mathbf Z$ for simplicity)
 and $\widetilde{\mathcal F}_{T-}:=\sigma\left(\bigcup_{t\in[0,T)}\widetilde\cF_t\right)$.
Moreover, by Proposition 1 in Fitzsimmons, Pitman and Yor \cite{FitPitYor},
 $(\widetilde\PP_{\mathbf a,\mathbf b}^T)_{\mathbf b\in\RR^d}$ is a regular version of the family of
 conditional distributions $\widetilde\PP(\cdot\mid {\mathbf Z}_T =\mathbf b)$, $\mathbf b\in\RR^d$,
 where \ $\widetilde\PP$ \ denotes the law of $(\mathbf Z_t)_{t\geq 0}$.
 Hence for any $n\in\NN$, $0<t_1<\cdots t_n<T$ and any $A\in\mathcal B(\RR^{nd})$,
 by Theorem \ref{gen_multivariate_intrep_bridge}, we get
 $$
  \widetilde\PP\big( ({\mathbf Z}_{t_1}^\top,\ldots,
                           {\mathbf Z}_{t_n}^\top)^\top \in A
                           \mid {\mathbf Z}_T = \mathbf b\big)
      =\widetilde\PP_{\mathbf a,\mathbf b}^T
       \big( ({\mathbf Z}_{t_1}^\top,\ldots,{\mathbf Z}_{t_n}^\top)^\top \in A\big)
   =\PP\big((\mathbf U_{t_1}^\top,\ldots,\mathbf U_{t_n}^\top)^\top\in A\big).
 $$

The next remark shows that the integral and anticipative representation of the
 bridge are quite different.

\begin{Rem}\label{Remark18}
Note that the process $(\mathbf Y_t)_{t\in[0,T]}$ defined in \eqref{anticipative_def} is only
 a weak solution of the SDE \eqref{gen_mult_bridge_sde}, since in contrast to the bridge
 $(\mathbf U_t)_{t\in[0,T]}$ it is not adapted to the filtration $(\mathcal F_t)_{t\geq0}$
 of the underlying Wiener process $\mathbf B$.
This can be easily seen by the definition of $\mathbf Y_t$
 which requires the knowledge of $\mathbf Z_T$ at any time point $t\in(0,T)$.
Nevertheless we have $\mathbf Y_t$ and $\mathbf Z_T$ are independent for any $t\in[0,T]$,
 since by part (a) of Lemma \ref{covariance_identities},
\begin{align*}
 \Cov(\mathbf Y_t,\mathbf Z_T)
 & =\Cov(\mathbf Z_t,\mathbf Z_T)-\Gamma(0,t)^\top\big(\Gamma(0,T)^\top\big)^{-1}\Cov(\mathbf Z_T,\mathbf Z_T)\\
 & =\Gamma(0,t)^\top E(T,0)^\top-\Gamma(0,t)^\top\big(\Gamma(0,T)^\top\big)^{-1}\Gamma(0,T)^\top E(T,0)^\top
  = 0\in\RR^{d\times d},
\end{align*}
and the random vector $(\mathbf Y_t^\top,\mathbf Z_T^\top)^\top$ has a Gauss distribution.
\end{Rem}

In the next remark we compare the SDE \eqref{gen_mult_bridge_sde} derived
 for the integral representation \eqref{gen_multivariate_intrep} of the bridge $\mathbf U$
 with the corresponding result of Delyon and Hu \cite{DelHu}.

\begin{Rem}\label{Remark14}
In this remark we discuss the connections between Proposition 3 in
Delyon and Hu \cite{DelHu} and our Theorem \ref{LEMMA_multi_DE}.
Let $p,d\in\NN$ and let us consider the SDE
\begin{align}\label{Delyon_Hu_egyenlet}
  \begin{cases}
   \dd \mathbf Z_t=(A(t)\mathbf Z_t + \mathbf g(t)+\sigma(t)\mathbf h(t,\mathbf Z_t))\,\dd t
                    + \sigma(t)\,\dd \mathbf B_t,\qquad t\geq0,\\
   \phantom{\dd} \mathbf Z_0=\mathbf a,
  \end{cases}
\end{align}
 where $\mathbf a\in\RR^d$, $A:\RR_+\to\RR^{d\times d},$ \ $\mathbf g:\RR_+\to\RR^d$ \ and
 \ $\sigma:\RR_+\to\RR^{d\times p}$ \ are continuous functions,
 \ $\mathbf h:\RR_+\times\RR^d\to\RR^p$ \ is a locally bounded function such that it is locally Lipschitz with
 respect to its second variable uniformly with respect to its first
 variable, i.e., for any \ $R>0,$ \ there exists a constant \ $C_R>0$ \
 such that for any \ $(t,\mathbf x), (t,\mathbf y)\in\RR_+\times\RR^d$ \ with \ $\Vert \mathbf x\Vert\leq R,$
 $\Vert \mathbf y\Vert\leq R$ \ we have
 $$
  \Vert \mathbf h(t,\mathbf x)-\mathbf h(t,\mathbf y)\Vert\leq C_R\Vert \mathbf x-\mathbf y\Vert,
   \quad \forall \;t\geq 0.
 $$
Moreover, we assume that the SDE \eqref{Delyon_Hu_egyenlet} has a (unique) strong solution
 (the question of uniqueness is not important here, Delyon and Hu \cite{DelHu} suppose
 only that there exists a strong solution), \ $\mathbf h$ \ is continuous with respect to its first variable,
 and that \ $\sigma$ \ admits a measurable left pseudo-inverse, denoted by
 \ $\sigma^+:=(\sigma^\top\sigma)^{-1}\sigma^\top$, \ which is left continuous
 (here for simplicity we suppose also that the symmetric matrix \ $\sigma^\top\sigma$ \ is
 positive definite).
Let us denote by \ $P_t,$ $t\geq 0,$ \ the unique solution of the deterministic matrix differential equation
 \ $P_t'=A(t)P_t,$ $t\geq 0,$ \ with initial condition \ $P_0=I_d$.
\ With the special choices \ $\mathbf h(t):= \mathbf 0$, $t\geq 0$,
 $\mathbf g(t):=\mathbf r(t),$ $t\geq 0,$
 $A(t):= Q(t),$ $t\geq 0,$ $\sigma(t):=\Sigma(t)$, $t\geq 0$, \ we get the
 SDE \eqref{Delyon_Hu_egyenlet} is the same as the SDE \eqref{gen_multivariate_system}
 with initial condition \ $\mathbf Z_0 = \mathbf a$.
\ Further, we have \ $P_t=E(t,0),$ $t\geq 0,$ \ and, by \eqref{evolution} and
 \eqref{evoinvert},
 \begin{align*}
   M_t& := \int_t^T P_u^{-1}\sigma(u)\sigma(u)^\top (P_u^{-1})^\top \,\dd u
         = \int_t^T E(u,0)^{-1}\Sigma(u)\Sigma(u)^\top (E(u,0)^{-1})^\top\,\dd u \\
       & = \int_t^T E(0,u) \Sigma(u)\Sigma(u)^\top E(0,u)^\top\,\dd u
         = E(0,t) \Gamma(t,T) E(0,T)^\top \\
       & = E(0,t)E(t,T) \kappa(t,T) E(0,T)^\top
         = E(0,T) \kappa(t,T) E(0,T)^\top,
         \qquad t\in[0,T].
 \end{align*}
Hence, by our assumption \eqref{kappa_assumption} on \ $\kappa$, \ $M_t$ \ is positive definite
 for all \ $t\in[0,T)$.
 Since \ $Q$, \ $\mathbf r$ \ and \ $\Sigma$ \ are continuous, if we suppose also that
 \ $\Sigma$ \ has a left continuous (measurable) left pseudo-inverse, which is guaranteed by assuming that
 $\Sigma\Sigma^\top$ is positive definite, then, by Proposition 3 in Delyon and Hu \cite{DelHu},
 the bridge \ $(\mathbf U_t)_{t\in[0,T]}$ \ from $\mathbf a$ to $\mathbf b$ over \ $[0,T]$ \
 derived from \ $\mathbf Z$ \ given by the SDE \eqref{gen_OU_egyenlet} (with initial condition
 $\mathbf Z_0=\mathbf a$) is a strong solution of the SDE
 \begin{align}\label{Delyon_Hu_egyenlet2}
  \begin{cases}
   \dd \mathbf U_t=A(t)\mathbf U_t\dd t + \mathbf g(t)\dd t + \sigma(t)\sigma(t)^\top (P_t^{-1})^\top M_t^{-1}
             \Big(P_t^{-1}(\EE \mathbf Z_t-\mathbf U_t)-P_T^{-1}(\EE \mathbf Z_T - \mathbf b)\Big)\,\dd t\\
   \phantom{\dd \mathbf U_t=\;}
             +\sigma(t)\,\dd \mathbf B_t,\qquad t\in[0,T),\\
   \phantom{\dd} \mathbf U_0 = \mathbf a.
  \end{cases}
 \end{align}
To be a little bit more precise, the definition of a bridge in Delyon and Hu \cite{DelHu}
 is different from our definition: they define a bridge as in Qian and Zheng \cite{QiaZhe}, Lyons and
 Zheng \cite{LyoZhe}, i.e., via Radon-Nycodim derivatives (detailed below),
 and using their Theorem 2 and Proposition 3 the bridge process is not a strong solution
 of the SDE \eqref{Delyon_Hu_egyenlet2}, but equals in law to this strong solution.
 We also note that the results of Qian and Zheng \cite{QiaZhe} and
 Lyons and Zheng \cite{LyoZhe} are valid for time-homogeneous diffusions,
 while Delyon and Hu \cite{DelHu} consider time inhomogeneous diffusions.
Further, Qian and Zheng \cite{QiaZhe} refer to their Section 2.1 on conditional
 processes as a set of folklore facts for which they could not find a reference.
In what follows we briefly describe a possible {\it heuristic} approach for
 time inhomogeneous diffusions which is a counterpart of the approach presented in
 Qian and Zheng \cite[Section 2.1]{QiaZhe} for time homogeneous diffusions.
Let \ $\PP_{\mathbf a, \mathbf b}^T$ \ be the probability measure on
 \ $\cF_{T-}=\sigma\left(\bigcup_{t\in[0,T)}\cF_t\right)$ \ defined by
 \begin{align*}
    \PP_{\mathbf a, \mathbf b}^T(F)
      := \int_F \frac{p^{\mathbf Z}_{0,t}(\mathbf a,\mathbf Z_t)p^{\mathbf Z}_{t,T}(\mathbf Z_t,\mathbf b)}
                     {p^{\mathbf Z}_{0,T}(\mathbf a,\mathbf b)}
                \,\dd \PP,
         \qquad F\in\cF_t,\;\; t\in[0,T),
 \end{align*}
 i.e., the Radon--Nycodim derivative of \ $\PP_{\mathbf a, \mathbf b}^T$ \ with respect to
 \ $\PP_{\mathbf a}^{\mathbf Z_t}$, \ considering them as probability measures on \ $(\Omega,\cF_t)$, \ is given by
  \[
     \frac{\dd\PP_{\mathbf a, \mathbf b}^T}{\dd\PP_{\mathbf a}^{\mathbf Z_t}}
             = \frac{p^{\mathbf Z}_{t,T}(\mathbf Z_t,\mathbf b)}{p^{\mathbf Z}_{0,T}(\mathbf a,\mathbf b)},
  \]
 where \ $\PP_{\mathbf a}^{\mathbf Z_t}$ \ is given by
 \[
   \PP_{\mathbf a}^{\mathbf Z_t}(F):=\int_F p^{\mathbf Z}_{0,t}(\mathbf a,\mathbf Z_t)
                                                   \,\dd\PP,
     \qquad F\in\cF_t.
 \]
Using that \ $\mathbf Z$ \ is almost surely (left) continuous we get \ $\cF_{T-}=\cF_T$
 \ (see, e.g., Karatzas and Shreve \cite[Problem 2.7.6 and Corollary 2.7.8]{KarShr}),
 and hence \ $\PP_{\mathbf a, \mathbf b}^T$ \ is a probability measure also on \ $\cF_T$.
A possible generalization of Lemma 2.1 in Qian and Zheng \cite{QiaZhe} for time inhomogeneous diffusions
 sounds as follows: under the probability measure \ $\PP_{\mathbf a, \mathbf b}^T$, \ the process
 \ $(\mathbf Z)_{t\in[0,T]}$ \ has a.s. continuous sample paths and admits transition densities
 given in \eqref{gen_multivariate_bridge_densities}, i.e.,
 \[
    \PP_{\mathbf a, \mathbf b}^T(Z_t\in B\mid Z_s) = \int_B p^{\mathbf U}_{s,t}(\mathbf Z_s,\mathbf y)\,\dd y,
    \qquad B\in\cB(\RR), \;\; 0\leq s<t<T,
 \]
 where \ $\cB(\RR)$ \ denotes the set of Borel sets in \ $\RR$.
\ Moreover, by the introduction in Delyon and Hu \cite{DelHu}, \ $\PP_{\mathbf a, \mathbf b}^T$ \ equals
 the law of the weak solution of the SDE
 \begin{align*}
   \dd {\mathbf V}_t
    = \Big( Q(t){\mathbf V}_t + {\mathbf r}(t) + \sigma(t)\sigma(t)^\top {\mathbf L}(t,{\mathbf V}_t)\Big)\dd t
      +\Sigma(t)\,\dd {\mathbf B}_t,\qquad t\in[0,T),
 \end{align*}
 where
  \begin{align*}
      {\mathbf L}(t,{\mathbf x})
          := \nabla_{\mathbf x}\ln p^{\mathbf Z}_{t,T}({\mathbf x},{\mathbf b}),
          \qquad t\in[0,T),\;\; {\mathbf x}\in\RR^d,
  \end{align*}
 with the notation
  \[
   \nabla_{\mathbf x}f({\mathbf x})
     :=\left(\frac{\partial f}{\partial x_1}({\mathbf x}),
                   \ldots,\frac{\partial f}{\partial x_d}({\mathbf x})\right)^\top,
     \qquad {\mathbf x}=(x_1,\ldots,x_d)\in\RR^d
  \]
  for a differentiable function \ $f:\RR^d\to\RR$.
\ Note that the above SDE is the SDE \eqref{gen_multivariate_system} with a modified drift function.
However, we can not address a rigorous proof of this approach and in fact this is one of the motivations
 for our different approach in the present paper.

The right hand side of the SDE \eqref{Delyon_Hu_egyenlet2} takes the form
 \begin{align*}
   Q(t)\mathbf U_t\dd t  + \mathbf r(t)\dd t + \Sigma(t)\,\dd \mathbf B_t
                       & + \Sigma(t)\Sigma(t)^\top (E(t,0)^{-1})^\top
                            (E(0,T)^\top)^{-1} \kappa(t,T)^{-1} E(0,T)^{-1} \\
      &\phantom{+\;} \times \Big(E(t,0)^{-1}(\mathbf m_{\mathbf a}(0,t) -\mathbf U_t)
                - E(T,0)^{-1}(\mathbf m_{\mathbf a}(0,T) -\mathbf b)\Big)\,\dd t.
 \end{align*}
Here the coefficient of \ $\dd t$ \ is equal to
 \begin{align*}
  &Q(t)\mathbf U_t + \mathbf r(t) + \Sigma(t)\Sigma(t)^\top E(0,t)^\top E(T,0)^\top \kappa(t,T)^{-1}E(T,0) \\
  &\phantom{Q(t)\mathbf U_t + \mathbf r(t) +\;}
      \times \Bigg[
         E(0,t)\left(E(t,0)\mathbf a + \int_0^t E(t,u)q(u)\,\dd u - \mathbf U_t\right) \\
  &\phantom{Q(t)\mathbf U_t + \mathbf r(t) + \Bigg[ \;\;\;}       
            - E(0,T)\left(E(T,0)\mathbf a + \int_0^T E(T,u)q(u)\,\dd u - \mathbf b\right)
       \Bigg] \\
  & = \Sigma(t)\Sigma(t)^\top E(0,t)^\top E(T,0)^\top \kappa(t,T)^{-1}E(T,0)\\
  &\phantom{=}\times\left(
              \int_0^tE(0,u)q(u)\,\dd u - E(0,t)\mathbf U_t
              - \int_0^T E(0,u)q(u)\,\dd u + E(0,T)\mathbf b
      \right) + Q(t)\mathbf U_t + \mathbf r(t) \\
  & = - \Sigma(t)\Sigma(t)^\top  E(T,t)^\top \kappa(t,T)^{-1} E(T,t)\mathbf U_t\\
  &\phantom{=\;}   - \Sigma(t)\Sigma(t)^\top E(T,t)^\top \kappa(t,T)^{-1} E(T,0) \int_t^T E(0,u)q(u)\,\dd u \\
  &\phantom{=\;}  + \Sigma(t)\Sigma(t)^\top E(T,t)^\top \kappa(t,T)^{-1}\mathbf b
     + Q(t)\mathbf U_t + \mathbf r(t) \\
  & =  - \Sigma(t)\Sigma(t)^\top  E(T,t)^\top \Gamma(t,T)^{-1}\mathbf U_t
     + \Sigma(t)\Sigma(t)^\top (\Gamma(t,T)^\top)^{-1} \mathbf m_{\mathbf b}^-(t,T)
     + Q(t)\mathbf U_t + \mathbf r(t).
 \end{align*}
This implies that the SDE \eqref{Delyon_Hu_egyenlet2} with our
special choices is the same as the SDE \eqref{gen_mult_bridge_sde}
 in Theorem \ref{LEMMA_multi_DE}.
But we emphasize again that the definition of a bridge in Delyon and Hu \cite{DelHu} is different
 from our definition and also the proofs of Proposition 3 in Delyon and Hu \cite{DelHu}
 and our Theorem \ref{LEMMA_multi_DE} are different.
Hence our Theorem \ref{LEMMA_multi_DE} is not an immediate consequence of Proposition 3
 in Delyon and Hu \cite{DelHu}.
\end{Rem}

In the next remark we compare the anticipative representation \eqref{anticipative_def} of the bridge
 \ $\mathbf U$ \ with the corresponding result in Delyon and Hu \cite{DelHu}.

\begin{Rem}\label{Remark12}
In this remark we discuss the connections between Theorem 2 in
  Delyon and Hu \cite{DelHu} and our Theorem \ref{anticipative_bridge}.
Let us consider the SDE \eqref{Delyon_Hu_egyenlet} in Remark
 \ref{Remark14}. With the special choices \ $\mathbf h(t):=\mathbf 0$, $t\geq 0$,
 $\mathbf g(t):=\mathbf r(t),$ $t\geq 0$, $A(t):=Q(t),$
 $t\geq 0$, $\sigma(t):=\Sigma(t)$, $t\geq 0,$ \ we get the SDE \eqref{Delyon_Hu_egyenlet}
 is the same as the SDE \eqref{gen_multivariate_system} with initial condition \ $\mathbf Z_0=\mathbf a$.
Theorem 2 in Delyon and Hu \cite{DelHu} implies that the linear process bridge from $\mathbf a$ to
 $\mathbf b$ over \ $[0,T]$ \ derived from \ $\mathbf Z^*$, which is given by the SDE
 \eqref{gen_multivariate_system} with initial condition \ $\mathbf Z_0^*=\mathbf a$, equals in law
 the process
 $$
   \mathbf Y_t^*
       = R^*(t,T)R^*(T,T)^{-1}\mathbf b
         +\Big(\mathbf Z_t^*- R^*(t,T) R^*(T,T)^{-1}\mathbf Z_T^*\Big),
           \quad t\in[0,T],
 $$
where \ $R^*(s,t):=\Cov(\mathbf Z_s^*,\mathbf Z_t^*),$
 $s,t\geq 0,$ \ is the covariance function of \ $\mathbf Z^*.$ \
Since \ $\mathbf Z_t^*=E(t,0)\mathbf a + \mathbf Z_t$, $t\geq 0$, \ where \ $(\mathbf Z_t)_{t\geq 0}$ \ is
 given by the SDE \eqref{gen_multivariate_system} with initial condition \ $\mathbf Z_0=\mathbf 0$, \
 by Lemma \ref{covariance_identities}, we get for all \ $t\in[0,T]$,
\begin{align*}
  \mathbf Y_t^*
      & = (E(T,0)\Gamma(0,t))^\top \big((E(T,0)\Gamma(0,T))^\top\big)^{-1}\mathbf b \\
      &\phantom{=\;}
         + \Big[ E(t,0)\mathbf a + \mathbf Z_t
               - (E(T,0)\Gamma(0,t))^\top \big((E(T,0)\Gamma(0,T))^\top\big)^{-1}
                  (E(T,0)\mathbf a+\mathbf Z_T)
            \Big]\\
      & = \Gamma(0,t)^\top (\Gamma(0,T)^\top)^{-1}\mathbf b
         + \Big(E(t,0) - \Gamma(0,t)^\top (\Gamma(0,T)^\top)^{-1} E(T,0)\Big)\mathbf a\\
      &\phantom{=\;}
         + \mathbf Z_t - \Gamma(0,t)^\top (\Gamma(0,T)^\top)^{-1} \mathbf Z_T \\
      & = \Gamma(t,T)\Gamma(0,T)^{-1}\mathbf a + \mathbf Z_t
           - \Gamma(0,t)^\top (\Gamma(0,T)^\top)^{-1}(\mathbf Z_T -\mathbf b)
         =Y_t,
\end{align*}
 where the last but one equality follows by part (c) of Lemma \ref{gen_sigmaintidentity}.
The proof of our Theorem \ref{anticipative_bridge} is the
 very same as the corresponding part of the proof of Theorem 2 in
 Delyon and Hu \cite{DelHu}.
\end{Rem}

\begin{Rem}\label{Remark17}
With the notations of Remark \ref{Remark12} one may define a bridge from $\mathbf 0$ to
 $\mathbf 0$ over \ $[0,T]$ \ derived from \ $\mathbf Z^*$ \ by
 \ $\mathbf Z_t^*-\EE(\mathbf Z_t^*\mid\mathbf Z_T^*)$, $t\in[0,T]$, \ subtracting
 from $\mathbf Z^*$ its conditional expectation given the process at the endpoint of the bridge.
Then, by Theorem 2 in Chapter II, \S13 of Shiryaev \cite{Shi} and our assumption
 \eqref{kappa_assumption}, it is known that the conditional distribution of \ $\mathbf Z_t^*$ \ given
 \ $\mathbf Z_T^*=\mathbf x$ \ is normal with mean
 \ $\EE\mathbf Z_t^*+ R^*(t,T) R^*(T,T)^{-1}(x-\EE\mathbf Z_T^*)$, $t\in[0,T].$
Hence we have
 \begin{align*}
  \mathbf Z_t^*-\EE(\mathbf Z_t^*\mid\mathbf Z_T^*)
    & =\mathbf Z_t^*-\EE\mathbf Z_t^*-R^*(t,T) R^*(T,T)^{-1}(\mathbf Z_T^*-\EE\mathbf Z_T^*)\\
    & =\mathbf Z_t^*- R^*(t,T) R^*(T,T)^{-1}\mathbf Z_T^*
       -\EE(\mathbf Z_t^*- R^*(t,T) R^*(T,T)^{-1}\mathbf Z_T^*)\\
    & =\mathbf Y_t^*-\EE\mathbf Y_t^*
        = \mathbf Y_t - \EE\mathbf Y_t,\qquad t\in[0,T],
\end{align*}
 which is nothing else but the centered anticipative representation
 of the bridge from $\mathbf 0$ to $\mathbf 0$.
Thus in general this definition of the bridge has a different mean function
 than the bridge given by Definition \ref{DEF_bridge_multidim}.
\end{Rem}

In case of dimension 1, we will also study the connections between Proposition 4
 in Gasbarra, Sottinen and Valkeila \cite{GasSotVal} and our Theorem \ref{anticipative_bridge}.
The reason for restricting ourselves to the case dimension one is that Gasbarra, Sottinen
 and Valkeila \cite{GasSotVal} consider only one-dimensional processes.

\section{One-dimensional linear process bridges}\label{SECTION_ONEDIMENSIONAL}

Let us consider a general one-dimensional linear process given by the linear SDE
 \begin{align}\label{gen_OU_egyenlet}
   \dd Z_t=\big(q(t)\,Z_t+r(t)\big)\,\dd t+\sigma(t)\,\dd B_t,\qquad t\geq0,
 \end{align}
 with continuous functions $q:\RR_{+}\to\RR,$  $\sigma:\RR_+\to\RR$ and
 $r:\RR_+\to\RR$, where $(B_t)_{t\geq 0}$ is a standard Wiener process.
By Section 5.6 in Karatzas and Shreve \cite{KarShr}, it is known that
 there exists a strong solution of the SDE \eqref{gen_OU_egyenlet}, namely
\begin{align}\label{gen_OU_egyenlet_megoldas}
  Z_t=\ee^{\bar q(t)}\left(Z_0+\int_0^t\ee^{-\bar q(s)}r(s)\,\dd s+\int_0^t\ee^{-\bar q(s)}\sigma(s)\,\dd B_s
                     \right),  \qquad t\geq 0,
 \end{align}
 with $\bar q(t):=\int_{0}^tq(u)\,\dd u$, $t\geq 0$, and strong uniqueness for the SDE
 \eqref{gen_OU_egyenlet} holds.
In what follows, we assume that $Z_0$ has a Gauss distribution independent of $(B_t)_{t\geq 0}$.
We call the process $(Z_t)_{t\geq 0}$ a one-dimensional linear process.
One can easily derive that for \ $0\leq s<t$ \ we have
 \begin{align}\label{SEGED26}
   Z_{t}=\ee^{\bar q(t)-\bar q(s)}Z_{s}+\int_s^t\ee^{\bar q(t)-\bar q(u)}r(u)\,\dd u
         +\int_{s}^t\ee^{\bar q(t)-\bar q(u)}\sigma(u)\,\dd B_u.
 \end{align}
 Hence, given \ $Z_s=x,$ \ the distribution of \ $Z_t$ \ does not depend on \ $(Z_{r})_{r\in[0,s)}$ \
 which yields that \ $(Z_{t})_{t\geq0}$ \ is a Markov process.
Moreover, for any \ $x\in\RR$ \ and \ $0\leq s<t$ \ the conditional distribution of \ $Z_t$ \ given
 \ $Z_s=x$ \ is Gauss with mean
 $$
   m_x(s,t):=\ee^{\bar q(t)-\bar q(s)}x+\int_s^t\ee^{\bar q(t)-\bar q(u)}r(u)\,\dd u,
 $$
 and with variance
 $$
   \gamma(s,t):=\int_{s}^t\ee^{2(\bar q(t)-\bar q(u))}\sigma^2(u)\,\dd u<\infty.
 $$
In what follows we put the following assumption
 \begin{align}\label{sigma_assumption}
   \sigma(t)\ne0\qquad \text{for all \ $t\geq 0$.}
 \end{align}
This yields that the variance \ $\gamma(s,t)$ \ is positive for all \ $0\leq s<t$
(which corresponds to condition \eqref{kappa_assumption} in dimension one).
Hence \ $(Z_{t})_{t\geq0}$ \ is a Gauss-Markov process with transition densities
  \begin{equation}\label{gen_OU_densities}
  p_{s,t}^Z(x,y)
  =\frac{1}{\sqrt{2\pi\gamma(s,t)}}
     \exp\left\{-\frac{(y-m_x(s,t))^2}{2\gamma(s,t)}\right\},
       \qquad 0\leq s<t,\quad x,y\in\RR.
 \end{equation}

For all \ $a,b\in\RR$ \ and \ $0\leq s\leq t<T$, \ let us introduce the notations
\begin{equation}\label{gen_OU_exp}
 n_{a,b}(s,t)
  :=\frac{\gamma(s,t)}{\gamma(s,T)}\,\ee^{\bar q(T)-\bar q(t)}
     \left(b-\int_t^T\ee^{\bar q(T)-\bar q(u)}r(u)\,\dd u\right)
     +\frac{\gamma(t,T)}{\gamma(s,T)}m_{a}(s,t),
 \end{equation}
 and
 \begin{equation}\label{gen_OU_var}
  \sigma(s,t):=\frac{\gamma(s,t)\,\gamma(t,T)}{\gamma(s,T)}.
 \end{equation}

Theorem \ref{gen_multivariate_intrep_bridge} has the following consequence.

\begin{Thm}\label{THEOREM2}
Let us suppose that condition \eqref{sigma_assumption} holds.
For fixed $a,b\in\RR$ and $T>0$, let the process $(U_t)_{t\in[0,T)}$ be given by
 \begin{align}\label{gen_OU_integral}
    U_t:= n_{a,b}(0,t)
                + \int_{0}^t\frac{\gamma(t,T)}{\gamma(s,T)}\ee^{\bar q(t)-\bar q(s)}\sigma(s)\,\dd B_{s},
                \quad t\in[0,T).
 \end{align}
Then for any $t\in[0,T)$ the distribution of $U_t$ is Gauss
 with mean $n_{a,b}(0,t)$ and with variance $\sigma(0,t).$
\ Especially, $U_t\to b$ almost surely (and hence in probability) and in $L^2$ as $t\uparrow T$.
Hence the process $(U_t)_{t\in[0,T)}$ can be extended to an almost surely (and hence stochastically) and
 $L^2$-continuous process $(U_t)_{t\in[0,T]}$ with $U_0=a $ and $U_T=b$.
Moreover, $(U_t)_{t\in[0,T]}$ is a Gauss-Markov process
 and for any $x\in\RR$ and $0\leq s<t<T$ the transition density
 $\RR\ni y\mapsto p_{s,t}^U(x,y)$ of $U_t$ given $U_s=x$ is given by
 \begin{align*}
   p_{s,t}^U(x,y)
     = \frac{1}{\sqrt{2\pi\sigma(s,t)}}
    \exp\left\{-\frac{\left(y-n_{x,b}(s,t)\right)^2}
                      {2\sigma(s,t)}\right\},
        \qquad y\in\RR.
 \end{align*}
 \end{Thm}

The proof of Theorem \ref{THEOREM2} can be found in the Appendix.

For completeness we formulate the definition of a one-dimensional linear process bridge,
 which definition is a special case of the multidimensional one  (see Definition \ref{DEF_bridge_multidim}).

\begin{Def}\label{DEF_bridge_onedim}
Let $(Z_t)_{t\geq0}$ be the one-dimensional linear process given by the SDE
 \eqref{gen_OU_egyenlet} with an initial Gauss random variable $Z_0$ independent
 of $(B_t)_{t\geq0}$ and let us assume that condition \eqref{sigma_assumption} holds.
For fixed $a, b\in\RR$ and $T>0$, the process $(U_{t})_{t\in[0,T]}$
 defined in Theorem \ref{THEOREM2} is called a linear process bridge
 from $a$ to $b$ over $[0,T]$ derived from \ $Z$.
 More generally, we call any almost surely continuous (Gauss) process on the time
 interval \ $[0,T]$ \ having the same finite-dimensional distributions as $(U_{t})_{t\in[0,T]}$
 a linear process bridge from $a$ to $b$ over $[0,T]$ derived from \ $Z$.
\end{Def}

Theorem \ref{LEMMA_multi_DE} has the following consequence.

\begin{Thm}\label{LEMMA4}
Let us suppose that condition \eqref{sigma_assumption} holds.
The process \ $(U_t)_{t\in[0,T)}$ \ defined by
\eqref{gen_OU_integral} is a unique strong solution of the linear SDE
 \begin{align}\label{gen_OU_hid1_egyenlet}
   \begin{split}
   \dd U_t= & \left(q(t)-\frac{\ee^{2(\bar q(T)-\bar q(t))}}{\gamma(t,T)}\sigma^2(t)\right)\,U_{t}\,\dd t\\
   & +\left(r(t)+\frac{\ee^{\bar q(T)-\bar q(t)}}{\gamma(t,T)}
       \left(b-\int_t^T\ee^{\bar q(T)-\bar q(u)}r(u)\,\dd u\right)\sigma^2(t)\right)\dd t
             +\sigma(t)\,\dd B_t
   \end{split}
 \end{align}
 for \ $t\in[0,T)$ \ and with initial condition \ $U_0=a$, \ and strong uniqueness
 for the SDE \eqref{gen_OU_hid1_egyenlet} holds.
\end{Thm}

As a consequence of Theorem \ref{anticipative_bridge} we give an anticipative representation
 of the linear process bridge introduced in Theorem \ref{THEOREM2} and Definition \ref{DEF_bridge_onedim}.

\begin{Thm}\label{PROPOSITION6}
Let \ $(Z_t)_{t\geq 0}$ \ be a linear process
 given by the SDE \eqref{gen_OU_egyenlet} with initial condition \ $Z_0=0$
 and let us suppose that condition \eqref{sigma_assumption} holds.
\ Then the process $(Y_t)_{t\in[0,T]}$ given by
 $$
   Y_t:=a\frac{\widetilde R(t,T)}{\widetilde R(0,T)}
          + Z_t
               - \frac{\widetilde R(0,t)}{\widetilde R(0,T)}(Z_T-b),
        \quad t\in[0,T],
 $$
equals in law the linear process bridge \ $(U_t)_{t\in[0,T]}$ \ from $a$ to $b$ over
 $[0,T]$ derived from the process $Z,$ \ where
$$
 \widetilde R(s,t):=\gamma(s,t)\ee^{\bar q(s)-\bar q(t)},
    \qquad 0\leq s\leq t\leq T.
$$
Moreover,
$$
  \widetilde R(s,t)
     =\ee^{\bar q(s)-\bar q(t)}R(t,t)-\ee^{\bar q(t)-\bar q(s)}R(s,s),
   \qquad 0\leq s\leq t\leq T,
$$
where \ $R$ \ denotes the covariance function of \ $Z,$ \ and
\begin{align}\label{SEGED16}
  \begin{split}
    Y_t&=a\left(\ee^{\bar q(t)}-\ee^{2\bar q(T)-\bar q(t)}\frac{\gamma(0,t)}{\gamma(0,T)}\right)
         +b\ee^{\bar q(T)-\bar q(t)}\frac{\gamma(0,t)}{\gamma(0,T)}\\
       &\phantom{=\;\;}
          +\left(Z_t-\ee^{\bar q(T)-\bar q(t)}\frac{\gamma(0,t)}{\gamma(0,T)}Z_T\right),
         \qquad t\in[0,T].
   \end{split}
\end{align}
\end{Thm}

We remark that the process \ $(Y_t)_{t\in[0,T]}$ \ in Theorem \ref{PROPOSITION6} can be written
 also in the form
\begin{align}\label{SEGED27}
  Y_t=a\left(\ee^{\bar q(t)}-\ee^{\bar q(T)}\frac{R(t,T)}{R(T,T)}\right)
      +b\frac{R(t,T)}{R(T,T)}
      +\left(Z_t-\frac{R(t,T)}{R(T,T)}Z_T\right),
      \quad t\in[0,T].
\end{align}
We also note that in Remark \ref{Remark7} we will give an illuminating explanation
 for the representation \eqref{SEGED16}.

As a consequence of Proposition \ref{bridge_conditioning} now we present a usual conditioning
 property for one-dimensional linear processes.

\begin{Pro}\label{PROPOSITION5}
Let \ $a,b\in\RR$ \ and \ $T>0$ \ be fixed. Let \ $(Z_t)_{t\geq 0}$
\ be the one-dimensional linear process given by the SDE
\eqref{gen_OU_egyenlet} with initial condition \ $Z_0=a$ \
and let us assume that condition \eqref{sigma_assumption} holds.
Let \ $n\in\NN$ \ and \ $0<t_1<t_2<\ldots<t_n<T.$ \ Then the conditional
distribution of \ $(Z_{t_1},\ldots,Z_{t_n})$ \ given \ $Z_T=b$ \
equals the distribution of \ $(U_{t_1},\ldots,U_{t_n}),$ \ where
 \ $(U_t)_{t\in[0,T]}$ \ is the linear process bridge from $a$ to $b$ over
  $[0,T]$ derived from $(Z_t)_{t\geq 0}.$
\end{Pro}

Next we give an illuminating explanation for the representation \eqref{SEGED16}
 in Theorem \ref{PROPOSITION6} (see Remark \ref{Remark7}), but preparatory
 we present a generalization of Lamperti transformation (see, e.g., Karlin and Taylor
 \cite[page 218]{KarTay2}) for one-dimensional linear processes.
This generalization may be known, but we were not able to find any reference,
 its proof can be found in the Appendix.

\begin{Pro}\label{LEMMA5}
 Let \ $(B_t^*)_{t\geq 0}$ \ be a standard Wiener process starting
 from $0$ and
 \begin{align*}
   Z_t^*:=
          m_0(0,t)+\ee^{\bar q(t)}B^*(\ee^{-2\bar q(t)}\gamma(0,t)),
            \quad t\geq 0.
 \end{align*}
 Then \ $(Z_t^*)_{t\geq 0}$ \ is a weak solution of the SDE
 \eqref{gen_OU_egyenlet} with initial condition \ $Z_0^*=0$.
\end{Pro}

\begin{Rem}\label{Remark7}
 Using Proposition \ref{LEMMA5} one can give an illuminating explanation
 for the representation \eqref{SEGED16} in Theorem
 \ref{PROPOSITION6}. By Problem 5.6.14 in Karatzas and Shreve
 \cite{KarShr}, the process \ $(\widehat U_t)_{t\in[0,T]}$ \ defined
 by
 $$
  \widehat U_t:=a\frac{T-t}{T}+b\frac{t}{T}
                +\left(\widehat B_t-\frac{t}{T}\widehat B_T\right),
                \quad t\in[0,T],
 $$
equals in law the Wiener bridge from $a$ to $b$ over \ $[0,T],$ \
where \ $(\widehat B_t)_{t\geq 0}$ \ is a standard Wiener process.
 Motivated by Lemma 1 in Papie\.{z} and Sandison \cite{PapSan} and
 Proposition \ref{LEMMA5}, first we will do the time
change \ $[0,T]\ni t\mapsto \ee^{-2\bar q(t)}\gamma(0,t),$ \ the rescaling
with coefficient \ $\ee^{\bar q(t)},$ \ and then the translation
with \ $m_0(0,t)$ \ for the process \ $(\widehat U_t)_{t\in[0,T]}.$
\ Namely, we consider the process
 \begin{align*}
  U_t^*&:=m_0(0,t)
         +\ee^{\bar q(t)}
           \left(a\frac{\ee^{-2\bar q(T)}\gamma(0,T)-\ee^{-2\bar q(t)}\gamma(0,t)}{\ee^{-2\bar q(T)}\gamma(0,T)}
                 +b\frac{\ee^{-2\bar q(t)}\gamma(0,t)}{\ee^{-2\bar q(T)}\gamma(0,T)}\right.\\
       &\phantom{:=m_0(0,t)+\ee^{\bar q(t)}\left(\;\;\right.}
        \left.+ \widehat B(\ee^{-2\bar q(t)}\gamma(0,t))
                        -\frac{\ee^{-2\bar q(t)}\gamma(0,t)}{\ee^{-2\bar q(T)}\gamma(0,T)}
                         \widehat B(\ee^{-2\bar q(T)}\gamma(0,T))
           \right),\quad t\in[0,T].
 \end{align*}
Then for all \ $t\in[0,T]$ \ we have
\begin{align*}
  U_t^*&=m_0(0,t)
        +a\left(\ee^{\bar q(t)}-\ee^{2\bar q(T)-\bar q(t)}\frac{\gamma(0,t)}{\gamma(0,T)}\right)
        +b\ee^{\bar q(T)}\ee^{\bar q(T)-\bar q(t)}\frac{\gamma(0,t)}{\gamma(0,T)}\\
       &\phantom{=\;\;}
         + \ee^{\bar q(t)}\widehat B(\ee^{-2\bar q(t)}\gamma(0,t))
                        -\ee^{\bar q(T)-\bar q(t)}\frac{\gamma(0,t)}{\gamma(0,T)}
                         \ee^{\bar q(T)}\widehat B(\ee^{-2\bar q(T)}\gamma(0,T))=
\end{align*}                         
\begin{align*}                   
      &=a\left(\ee^{\bar q(t)}-\ee^{2\bar q(T)-\bar q(t)}\frac{\gamma(0,t)}{\gamma(0,T)}\right)
        +\big(\ee^{\bar q(T)}b+m_0(0,T)\big)\ee^{\bar q(T)-\bar q(t)}\frac{\gamma(0,t)}{\gamma(0,T)}\\
      &\phantom{=\;\;}
       +\Big(Z_t^*-\ee^{\bar q(T)-\bar q(t)}\frac{\gamma(0,t)}{\gamma(0,T)}Z_T^*\Big),
\end{align*}
where, using Proposition \ref{LEMMA5}, \ $(Z_t^*)_{t\geq 0}$ \ equals in
 law the one-dimensional linear process given by the SDE
 \eqref{gen_OU_egyenlet} with initial condition \ $Z_0=0$.
By Theorem \ref{PROPOSITION6}, the process \ $(U_t^*)_{t\in[0,T]}$ \ equals in law
 the one-dimensional linear process bridge \ $(U_t)_{t\in[0,T]}$ \
 from $a$ to \ $\ee^{\bar q(T)}b+m_0(0,T)$ \ over \ $[0,T]$ \ derived from \
 $Z$ \ given by the SDE \eqref{gen_OU_egyenlet} with initial condition \ $Z_0=0$.
Roughly speaking, we have to apply the same time change,
rescaling and translation to the Wiener bridge from $a$ to $b$ over
\ $[0,T]$ \ in order to get the linear process bridge
 from $a$ to \ $\ee^{\bar q(T)}b+m_0(0,T)$ \ over \ $[0,T]$ \
 (derived from $Z$ given by the SDE \eqref{gen_OU_egyenlet} with initial condition \ $Z_0=0$)
 what we apply to a Wiener process in order to get the linear process \ $Z$.

Especially, concerning Wiener bridges and Ornstein-Uhlenbeck bridges, we have to
 apply the same time change and rescaling to the Wiener bridge from $a$ to $b$ over \ $[0,T]$ \
in order to get the Ornstein-Uhlenbeck bridge from $a$ to
 \ $\ee^{qT}b$ \ over \ $[0,T]$ \ (derived from $Z$ given by the SDE
 \eqref{OU_egyenlet}) what we apply to a Wiener process in order
 to get the Ornstein-Uhlenbeck process \ $Z$.
 We note that the original definition of an Ornstein-Uhlenbeck bridge
 of Papie\.{z} and Sandison is different from ours, they define the
 bridge as a probability measure on the space of continuous
 functions \ $f:[0,T]\to\RR$ \ such that \ $f(0)=a$ \ and \ $f(T)=\ee^{qT}b.$
\end{Rem}

Next we formulate special cases of the presented one-dimensional results.

\begin{Rem}
Note that in case of $q(t)=q\ne 0$, $t\geq 0$, and $\sigma(t)=\sigma\ne 0,$ $t\geq 0$,
 (for any continuous deterministic forcing term $r$) the variance $\sigma(s,t)$
 defined by \eqref{gen_OU_var} gives back \eqref{OU_bridge_var}.
\proofend
\end{Rem}

Theorem \ref{THEOREM2} has the following consequence.

\begin{Rem}\label{REMARK_gen_transition}
Note that in case of \ $q(t)=q=0,$ \ $\sigma(t)=\sigma\not=0,$ \ $r(t)=0$, $t\geq 0$, \ and $a=0=b,$ \
 we recover the Wiener bridge \ $(\widetilde U_t)_{t\in[0,T]}$ \ from \ $0$ \
 to \ $0$ \ stated in \eqref{W_bridge}.
Moreover, in case of \ $q(t)=q\not=0,$ \ $\sigma(t)=\sigma\not= 0,$ \ and \ $r(t)=0$, $t\geq 0$,
 \ the linear process bridge (Ornstein-Uhlenbeck bridge) \ $(U_t)_{t\in[0,T]}$ \ from $a$ to $b$
 over \ $[0,T]$ \ defined in \eqref{gen_OU_integral} has the form
  \begin{align}\label{OU_hid1_intrep}
   U_t=a\,\frac{\sinh(q(T-t))}{\sinh(qT)}
        +b\,\frac{\sinh(qt)}{\sinh(qT)}
        +\sigma\int_0^t\frac{\sinh(q(T-t))}{\sinh(q(T-s))}\,\dd B_s,
         \qquad t\in[0,T),
 \end{align}
 and admits transition densities
   $$
    p_{s,t}^U(x,y)
       =\frac{1}{\sqrt{2\pi\sigma(s,t)}}
        \exp\left\{-\frac{\left(y-\frac{\sinh(q(t-s))}
                   {\sinh(q(T-s))}\,b-\frac{\sinh(q(T-t))}{\sinh(q(T-s))}\,x\right)^2}
                    {2\sigma(s,t)}\right\}
   $$
 for all \ $0\leq s<t<T$ \ and \ $x,y\in\RR,$ \ where \ $\sigma(s,t)$ \ is
 given by \eqref{OU_bridge_var}.
\end{Rem}

Theorem \ref{LEMMA4} has the following consequence.

\begin{Rem}\label{Remark_invariance}
Note that in case of \ $q(t)= q\not =0,$ \ $\sigma(t)=\sigma\not =0$ \ and \ $r(t)=0$, $t\geq 0$,
 \ the SDE \eqref{gen_OU_hid1_egyenlet} has the form
 \begin{align}\label{OU_hid1_egyenlet}
  \begin{cases}
   \dd U_t=q\left(-\coth(q(T-t))\,U_t+\frac{b}{\sinh(q(T-t))}\right)\,\dd t+\sigma\,\dd B_t,\qquad t\in[0,T),\\
   \phantom{\dd} U_0=a.
  \end{cases}
  \end{align}
  Note also that both the SDE \eqref{OU_hid1_egyenlet} and the integral
representation \eqref{OU_hid1_intrep} are invariant under a change
of sign for the parameter $q$. Hence the Ornstein-Uhlenbeck bridges
derived from the SDE \eqref{OU_egyenlet} with $q$ and $-q$ are (almost
surely) pathwise identical.
\end{Rem}

Theorem \ref{PROPOSITION6} has the following consequence.

\begin{Rem}\label{Remark4}
 We consider a special case of Theorem \ref{PROPOSITION6},
 namely, let us suppose that \ $r(t)=0$, $t\geq 0$,  and that there exist
 real numbers \ $q\not=0$ \ and \ $\sigma\not=0$ \ such that \
 $q(t)=q,$ $t\geq 0,$ \ and \ $\sigma(t)=\sigma,$ $t\geq 0.$ \
 Then \ $\bar q(t)=qt,$ $t\geq 0,$ \ and
\begin{align}\label{SEGED23}
 \begin{split}
  \widetilde R(s,t)
       &=\gamma(s,t)\ee^{\bar q(s)-\bar q(t)}
       =\sigma^2\ee^{q(s-t)}\int_s^t\ee^{2q(t-u)}\;\dd u
       =\sigma^2\ee^{q(s-t)}\frac{1}{2q}(\ee^{2q(t-s)}-1)\\
      &=\frac{\sigma^2}{2q}(\ee^{q(t-s)}-\ee^{-q(t-s)})
       =\frac{\sigma^2}{q}\sinh(q(t-s)),
      \quad 0\leq s\leq t\leq T,
  \end{split}
\end{align}
 and
 $$
     R(s,t)=\Cov(Z_s,Z_t)
           =\frac{\sigma^2}{2q}\ee^{q(s+t)}(1-\ee^{-2qs})
           =\frac{\sigma^2}{q}\ee^{qt}\sinh(qs),
       \quad 0\leq s\leq t.
 $$
An easy calculation shows that for all \ $t\in[0,T],$
 \begin{align*}
   &\frac{\widetilde R(0,t)}{\widetilde R(0,T)}
      =\ee^{q(T-t)}\frac{R(t,t)}{R(T,T)}
      =\frac{R(t,T)}{R(T,T)},\\[1mm]
   &\frac{\widetilde R(t,T)}{\widetilde R(0,T)}
      =\ee^{qt}-\ee^{2qT-qt}\frac{R(t,t)}{R(T,T)}
      =\frac{\sigma^2}{2q}\frac{\ee^{qt+2qT}(1-\ee^{-2qT})-\ee^{2qT+qt}(1-\ee^{-2qt})}{R(T,T)}\\
  &\phantom{\frac{\widetilde R(t,T)}{\widetilde R(0,T)}\;}
      =\frac{\sigma^2}{2q}\frac{\ee^{2qT-qt}-\ee^{qt}}{R(T,T)}
      =\frac{R(T-t,T)}{R(T,T)}.
\end{align*}
Hence the process \ $(Y_t)_{t\in[0,T]}$ \ introduced in Theorem
 \ref{PROPOSITION6} (with our special choices of \ $q,$ $r$ \ and \ $\sigma$) \
 has the form
$$
   Y_t=a\frac{R(T-t,T)}{R(T,T)}
         +b\frac{R(t,T)}{R(T,T)}
         +\left(Z_t-\frac{R(t,T)}{R(T,T)}Z_T\right),
        \quad t\in[0,T].
 $$
Moreover, by \eqref{SEGED23},
\begin{align}\label{SEGED16_old}
   Y_t
      =a\frac{\sinh(q(T-t))}{\sinh(qT)}
       +b\frac{\sinh(qt)}{\sinh(qT)}
       +\left(Z_t-\frac{\sinh(qt)}{\sinh(qT)}Z_T\right),
       \quad t\in[0,T].
\end{align}
 Finally, we remark that in case of \ $q(t)=q\ne 0,$ $\sigma(t)=\sigma\ne 0$, $t\geq 0$ \
 and \ $r(t)=0$, $t\geq 0$ \ with the special choices \ $q=-\sqrt{k\gamma}/2$ \
 and \ $\sigma=k/4,$ \ where \ $k>0$ \ and \ $\gamma>0,$ \ our Theorem \ref{PROPOSITION6}
 is the same as Lemma 1 in Papie\.{z} and Sandison \cite{PapSan}.
\end{Rem}

In the next remark we discuss the connections between Propositions 4 and 9
 in Gasbarra, Sottinen and Valkeila \cite{GasSotVal} and our Theorems \ref{PROPOSITION6}
 and \ref{THEOREM2} (anticipative and integral representation of the bridge in case of dimension one),
 respectively.

\begin{Rem}\label{Remark15}
It turns out that our Theorem \ref{PROPOSITION6} can be considered as a consequence of Proposition 4
 in Gasbarra, Sottinen and Valkeila \cite{GasSotVal}.
 Namely, by Theorem \ref{PROPOSITION6}, the process $(Y_t)_{t\in[0,T]}$ given by \eqref{SEGED27}
 equals in law the one-dimensional linear process bridge from \ $a$ \ to \ $b$ \
 over \ $[0,T]$ \ derived from the process \ $Z$ \ given by the SDE \eqref{gen_OU_egyenlet}
 with initial condition \ $Z_0=0$.
\ By Proposition 4 in Gasbarra, Sottinen and Valkeila \cite{GasSotVal}, if
\ $(Z_t^*)_{t\geq 0}$ \ is given by the SDE
\begin{align*}
  \begin{cases}
   \dd Z_t^*=(q(t)Z_t^*+r(t))\,\dd t +\sigma(t)\,\dd B_t,\qquad t\geq0,\\
   \phantom{\dd} Z_0^*=a,
  \end{cases}
\end{align*}
 then the bridge form \ $a$ \ to \ $b$ \ over \ $[0,T]$ \ derived from \ $Z^*$ \
 (defined as a Gauss process in the sense of Definition 2 in Gasbarra, Sottinen
 and Valkeila \cite{GasSotVal}) admits the representation
 $$
    b\frac{R^*(t,T)}{R^*(T,T)}+\left(Z_t^*-\frac{R^*(t,T)}{R^*(T,T)}Z_T^*\right),
       \qquad t\in[0,T],
 $$
  where \ $R^*(s,t):=\Cov(Z_s^*,Z_t^*,)$ $s,t\geq 0,$ \ is the covariance function
  of \ $Z^*.$
\ Since \ $Z^*_t=\ee^{\overline q(t)}a+Z_t,$ $t\geq 0,$ \ and
 \ $R^*(s,t)=R(s,t),$ $s,t\geq 0$ \ (where \ $R$ \ denotes the covariance function of \ $Z$),
 we have for all \ $t\in[0,T],$
 \begin{align*}
    b\frac{R^*(t,T)}{R^*(T,T)}&+Z_t^*-\frac{R^*(t,T)}{R^*(T,T)}Z_T^*\\
      &=a\left(\ee^{\overline q(t)}-\frac{R(t,T)}{R(T,T)}\ee^{\overline q(T)}\right)
       +b\frac{R(t,T)}{R(T,T)}
       +\left(Z_t-\frac{R(t,T)}{R(T,T)}Z_T\right),
 \end{align*}
as desired. But we emphasize that the definition of the bridge in Definition 2
in Gasbarra, Sottinen and Valkeila \cite{GasSotVal} is via a properly defined
conditional probability measure on the given probability space, and hence it differs
from our definition of a bridge. Hence in this respect our Theorem \ref{PROPOSITION6}
is not an immediate consequence of Proposition 4 in Gasbarra, Sottinen and Valkeila \cite{GasSotVal}.

 Further, if \ $(Z_t)_{t\geq 0}$ \ is given by the SDE \eqref{gen_OU_egyenlet}
 with \ $r(t)=0$, $t\geq 0$, \ and with initial condition \ $Z_0=0$, \ then the process
 \begin{equation}\label{contGmart}
   M_t:= a + \ee^{-\bar q(t)}Z_t = a + \int_0^t\ee^{-\bar q(s)}\sigma(s)\,\dd B_s,\quad t\geq0,
 \end{equation}
 is a continuous Gauss martingale (with respect to the filtration \ $(\cF_t)_{t\geq 0}$ \
 introduced in Section \ref{SECTION_MULTI}).
By Proposition 9 in Gasbarra, Sottinen and Valkeila \cite{GasSotVal}, the anticipative representation
 of the bridge from $a$ to $\ee^{-\bar q(T)}b$ over $[0,T]$ derived from \ $M$ \ is given by
\begin{align}\label{SEGED19}
  \widetilde Y_t := \ee^{-\bar q(T)}b\,\frac{\langle M\rangle_t}{\langle M\rangle_T}
      +\left(M_t-\frac{\langle M\rangle_t}{\langle M\rangle_T}\,M_T\right),
        \quad t\in[0,T],
\end{align}
 where \ $\langle M\rangle$ \ denotes the quadratic variation process of \ $M$.
Since
 \begin{equation}\label{contGmartquad}
  \langle M\rangle_t
     =\int_0^t\ee^{-2\bar q(s)}\sigma^2(s)\;\dd\langle B\rangle_s=\int_0^t\ee^{-2\bar q(s)}\sigma^2(s)
      \;\dd s=\ee^{-2\bar q(t)}\gamma(0,t), \qquad t\in[0,T],
 \end{equation}
\eqref{SEGED19} yields for $t\in[0,T],$
\begin{equation}\label{SEGED19a}
 \widetilde Y_t
  =\ee^{\bar q(T)-2\bar q(t)}b\,\frac{\gamma(0,t)}{\gamma(0,T)}+\ee^{-\bar q(t)}Z_t
   -\ee^{\bar q(T)-2\bar q(t)}\frac{\gamma(0,t)}{\gamma(0,T)}\,Z_T
   +a\left(1-\ee^{2(\bar q(T)-\bar q(t))}\frac{\gamma(0,t)}{\gamma(0,T)}\right).
\end{equation}
This implies that $(\ee^{\bar q(t)}\widetilde Y_t)_{t\in[0,T]}$ coincides with \eqref{SEGED16},
 the anticipative representation of the bridge from $a$ to $b$ over $[0,T]$ derived from $Z$.

Moreover, motivated by equation (3.2) in Proposition 9 of Gasbarra, Sottinen and Valkeila
 \cite{GasSotVal} (in case of \ $Z_0=0$ \ and \ $r(t)=0$, $t\geq 0$) a non-anticipative (integral) representation of the bridge from $a$ to $\ee^{-\bar q(T)}b$ over $[0,T]$ derived from $M$ is given by
\begin{align}\label{SEGED19b}
  \widetilde U_t:=a+\left(\ee^{-\bar q(T)}b-a\right)\,\frac{\langle M\rangle_t}{\langle M\rangle_T}
      +\int_0^t\frac{\langle M\rangle_{T,t}}{\langle M\rangle_{T,s}}\,\dd M_s,
        \quad t\in[0,T],
\end{align}
 where, similarly to \eqref{contGmartquad}, \ $\langle M\rangle_{T,t}$, $t\in[0,T]$, is given by
$$
  \langle M\rangle_{T,t}
     := \langle M\rangle_T - \langle M\rangle_t=
     \int_t^T\ee^{-2\bar q(s)}\sigma^2(s)\;\dd\langle B\rangle_s=\int_t^T\ee^{-2\bar q(s)}\sigma^2(s)\;\dd s=\ee^{-2\bar q(T)}\gamma(t,T).
$$
Hence \eqref{SEGED19b}, \eqref{contGmart} and \eqref{contGmartquad} yield for $t\in[0,T]$
\begin{equation*}
\widetilde U_t=a\left(1-\ee^{2(\bar q(T)-\bar q(t))}\frac{\gamma(0,t)}{\gamma(0,T)}\right)
           +\ee^{\bar q(T)-2\bar q(t)}b\,\frac{\gamma(0,t)}{\gamma(0,T)}
           +\int_0^t\frac{\gamma(t,T)}{\gamma(s,T)}\,\ee^{-\bar q(s)}\sigma(s)\;\dd B_s,
\end{equation*}
 which implies that \ $(\ee^{\bar q(t)}\widetilde U_t)_{t\in[0,T]}$ \ coincides with
 \eqref{gen_OU_integral}, the integral representation of the bridge from $a$ to $b$
 over $[0,T]$ derived from $Z$, since from the definition of $(\gamma(s,t))_{0\leq s\leq t}$
 one can easily derive that
 $$
  \ee^{\bar q(t)}-\ee^{2\bar q(T)-\bar q(t)}\frac{\gamma(0,t)}{\gamma(0,T)}
   =\ee^{\bar q(t)}\frac{\gamma(t,T)}{\gamma(0,T)},\qquad t\in[0,T).
 $$
\end{Rem}

\section{Appendix}\label{SECTION_PROOFS}

First we give sufficient conditions for positive definiteness of the Kalman matrices
 introduced in Section \ref{SECTION_MULTI}, see, e.g., Theorems 7.7.1-7.7.3 in Conti \cite{Con}.

\begin{Pro}\label{PROPOSITION_kappa}
Let \ $0\leq s<t$ be given. Then \ $\kappa(s,t)$ \
is positive definite if one of the following conditions is satisfied:
\begin{enumerate}
\item[(a)] there exists \ $t_{0}\in(s,t)$ \ such that \ $\Sigma(t_{0})$ \ has full rank $d$
           (for which $p\geq d$ is required).
\item[(b)] there exist \ $t_{0}\in(s,t),$ \ an open neighborhood \ $I_{0}$ \ around \ $t_{0}$ \
           and some $k\in\NN$ \ such that
           \ $\Sigma\in\mathcal C^{(k)}_{d\times p}(I_{0}),$ \ $Q\in\mathcal C^{(k-1)}_{d\times d}(I_{0})$ \
           and the controllability matrix
           \ $\big[\Sigma(t_{0}),\Delta\Sigma(t_{0}),\ldots,\Delta^k\Sigma(t_{0})\big]$ \
            has full rank $d$, where \ $\Delta$ \ is the operator
            \ $\Delta\Sigma(t)=\Sigma'(t)-Q(t)\Sigma(t),$ $t\in I_0$ \ and for all \ $n,m\in\NN$,
            \ $\mathcal C^{(k)}_{n\times m}(I_{0})$ \
            denotes the set of $k$-times continuously differentiable functions on \ $I_0$ \ with values
            in \ $\RR^{n\times m}$ (for which $(k+1)p\geq d$ is required).
\item[(c)]  there exist \ $t_{0}\in(s,t),$ \ an open neighborhood \ $I_{0}$ \ around \ $t_{0}$ \
            and some $k\in\NN$ \ such that
            \ $\Sigma\in\mathcal C^{(\infty)}_{d\times p}(I_{0}),$ \
            $Q\in\mathcal C^{(\infty)}_{d\times d}(I_{0})$ \
            and the controllability matrix
            \ $\big[\Sigma(t_{0}),\Delta\Sigma(t_{0}),\ldots,\Delta^k\Sigma(t_{0})\big]$ \
            has full rank $d$, \ where \ $\mathcal C^{(\infty)}_{n\times m}(I_{0})$ \
            denotes the set of infinitely differentiable functions on \ $I_0$ \ with values
            in \ $\RR^{n\times m}$ (for which $(k+1)p\geq d$ is required).
\end{enumerate}
\end{Pro}

Next we present two lemmata which will be used several times in the proofs later on.

\begin{Lem}\label{gen_sigma_identity}
Let us suppose that condition \eqref{kappa_assumption} holds.
For fixed $T>0$ and all $0\leq s<t<T$,
 \begin{equation}\label{gen_sigma_inverse}
  \Sigma(s,t)^{-1}=\kappa(s,t)^{-1}+E(T,t)^\top\kappa(t,T)^{-1}E(T,t).
 \end{equation}
Especially, $\Sigma(s,t)$ is symmetric and positive definite for all $0\leq s<t<T$.
\end{Lem}

\noindent{\bf Proof.}
By assumption \eqref{kappa_assumption},
 $\kappa(s,t)$ is symmetric and positive definite for all $0\leq s<t$, which implies
 that the right-hand side of \eqref{gen_sigma_inverse}, and
 thus also its inverse, is symmetric and positive definite. For $0\leq s<t<T$ we calculate
\begin{align*}
& \kappa(s,t)^{-1}+E(T,t)^\top\kappa(t,T)^{-1}E(T,t)\\
& \quad=\kappa(s,t)^{-1}\left(I_d+\kappa(s,t)E(T,t)^\top\kappa(t,T)^{-1}E(T,t)\right)\\
& \quad=\kappa(s,t)^{-1}\Bigg(I_d+\int_s^tE(t,u)\Sigma(u)\Sigma(u)^\top E(T,u)^\top\,{\rm d} u\\
& \qquad\qquad\qquad\qquad\left.\times\left(\int_t^TE(t,u)\Sigma(u)\Sigma(u)^\top E(T,u)^\top\,{\rm d} u\right)^{-1}\right)\\
& \quad=\kappa(s,t)^{-1}\int_s^TE(t,u)\Sigma(u)\Sigma(u)^\top E(T,u)^\top\,{\rm d} u\,\Gamma(t,T)^{-1}\\
&\quad=\Gamma(s,t)^{-1}\Gamma(s,T)\Gamma(t,T)^{-1}=\Sigma(s,t)^{-1},
\end{align*}
 which concludes the proof.
\proofend

\begin{Lem}\label{gen_sigmaintidentity}
Let us suppose that condition \eqref{kappa_assumption} holds.
For fixed $T>0$ and all $0\leq s<t<T$ we have
\begin{enumerate}
\item[(a)] $\kappa(t,T)^{-1}-\kappa(s,T)^{-1}=\Gamma(s,T)^{-1}\Gamma(s,t)\big(\Gamma(t,T)^\top\big)^{-1}$,
\item[(b)] $\Sigma(s,t)=\Gamma(t,T)\displaystyle\int_{s}^t\Gamma(u,T)^{-1}\Sigma(u)\Sigma(u)^\top\big(\Gamma(u,T)^\top\big)^{-1}\,{\rm d}u\,\Gamma(t,T)^\top$,
\item[(c)] $E(t,0)-\Gamma(0,t)^\top\big(\Gamma(0,T)^\top\big)^{-1}E(T,0)=\Gamma(t,T)\Gamma(0,T)^{-1}$,
\item[(d)] $\Gamma(s,T)^\top E(t,s)^\top-E(T,s)\Gamma(s,t)=\Gamma(t,T)^\top$.
\end{enumerate}
\end{Lem}
\noindent{\bf Proof.}
Since $\kappa(t,T)$ is symmetric we calculate
\begin{align*}
\kappa(t,T)^{-1}-\kappa(s,T)^{-1}& =\big(\kappa(t,T)^\top\big)^{-1}-\kappa(s,T)^{-1}\\
& =E(t,T)^\top\big(\Gamma(t,T)^\top\big)^{-1}-\Gamma(s,T)^{-1}E(s,T)\\
& =\Gamma(s,T)^{-1}\left(\Gamma(s,T)E(t,T)^\top-E(s,T)\Gamma(t,T)^\top\right)\big(\Gamma(t,T)^\top\big)^{-1},
\end{align*}
where the middle factor coincides with
$$\int_{s}^TE(s,u)\Sigma(u)\Sigma(u)^\top E(t,u)^\top\,{\rm d} u-\int_{t}^TE(s,u)\Sigma(u)\Sigma(u)^\top E(t,u)^\top\,{\rm d} u=\Gamma(s,t),$$
which proves (a). In a similar manner one can prove (c) and (d).
To prove (b), note that the function $(0,T)\ni u\mapsto\kappa(u,T)\in\RR^{d\times d}$ is a differentiable
curve in the set of symmetric and positive definite $(d\times d)$-matrices with
\begin{align*}
-\partial_1\kappa(u,T)
 & =-\frac{{\rm d}}{{\rm d}u}\int_{u}^T\!\!E(T,v)\Sigma(v)\Sigma(v)^\top E(T,v)^\top{\rm d}v
   =E(T,u)\Sigma(u)\Sigma(u)^\top E(T,u)^\top,\; u\in(0,T).
\end{align*}
Using (a) and the fact that
 $\partial_1\big(\kappa(u,T)^{-1}\big)=-\kappa(u,T)^{-1}\big(\partial_1\kappa(u,T)\big)\kappa(u,T)^{-1}$,
 $u\in(0,T)$, \ see, e.g., formula (3.2) on page 73 in Baker \cite{Bak}, we calculate
\begin{align*}
& \int_{s}^t\Gamma(u,T)^{-1}\Sigma(u)\Sigma(u)^\top\big(\Gamma(u,T)^\top\big)^{-1}\,{\rm d} u\\
& \quad=\int_{s}^t\kappa(u,T)^{-1}E(T,u)\Sigma(u)\Sigma(u)^\top E(T,u)^\top\kappa(u,T)^{-1}\,{\rm d} u\\
& \quad=-\int_{s}^{t}\kappa(u,T)^{-1}\big(\partial_1\kappa(u,T)\big)\kappa(u,T)^{-1}\,{\rm d} u\\
& \quad=\kappa(t,T)^{-1}-\kappa(s,T)^{-1}=\Gamma(s,T)^{-1}\Gamma(s,t)\big(\Gamma(t,T)^\top\big)^{-1},
\end{align*}
from which the assertion of part (b) easily follows.
\proofend

\vspace{1mm}

\noindent{\bf Proof of Lemma \ref{gen_multivariate_bridge}.}
 Since $\det\Sigma(s,t)=\det\kappa(s,t)\det\kappa(t,T)\big(\det\kappa(s,T)\big)^{-1}$, for $\mathbf x,\mathbf y\in\RR^d$ and $0\leq s<t<T$ we obtain, by \eqref{gen_multivariate_density},
$$\frac{p_{s,t}^{\mathbf Z}(\mathbf x,\mathbf y)\,p_{t,T}^{\mathbf Z}(\mathbf y,\mathbf b)}{p_{s,T}^{\mathbf Z}(\mathbf x,\mathbf b)}=\frac1{\sqrt{(2\pi)^d\det\Sigma(s,t)}}\,\exp\left\{-\frac12\,\psi_{s,t}(\mathbf x,\mathbf y)\right\},$$
where using \eqref{gen_sigma_inverse} and that
$$
  \mathbf b-\mathbf m_{\mathbf y}(t,T)=-E(T,t)\mathbf y+\mathbf m_{\mathbf b}^-(t,T),
$$
we calculate
\begin{align*}
\psi_{s,t}(\mathbf x,\mathbf y) & =\left\langle\kappa(s,t)^{-1}\big(\mathbf y-\mathbf m_{\mathbf x}(s,t)\big),\mathbf y-\mathbf m_{\mathbf x}(s,t)\right\rangle\\
& \quad+\left\langle\kappa(t,T)^{-1}\big(\mathbf b-\mathbf m_{\mathbf y}(t,T)\big),\mathbf b-\mathbf m_{\mathbf y}(t,T)\right\rangle\\
& \quad-\left\langle\kappa(s,T)^{-1}\big(\mathbf b-\mathbf m_{\mathbf x}(s,T)\big),\mathbf b-\mathbf m_{\mathbf x}(s,T)\right\rangle\\
& =\left\langle\left(\kappa(s,t)^{-1}+E(T,t)^\top\kappa(t,T)^{-1}E(T,t)\right)\mathbf y,\mathbf y\right\rangle\\
 & \quad+\left\langle\kappa(s,t)^{-1}\mathbf m_{\mathbf x}(s,t),\mathbf m_{\mathbf x}(s,t)\right\rangle
  -2\left\langle\kappa(s,t)^{-1}\mathbf y,\mathbf m_{\mathbf x}(s,t)\right\rangle\\
& \quad+\left\langle\kappa(t,T)^{-1}\mathbf m_{\mathbf b}^-(t,T),\mathbf m_{\mathbf b}^-(t,T)\right\rangle-2\left\langle\kappa(t,T)^{-1}E(T,t)\mathbf y,\mathbf m_{\mathbf b}^-(t,T)\right\rangle\\
& \quad-\left\langle\kappa(s,T)^{-1}\big(\mathbf b-\mathbf m_{\mathbf x}(s,T)\big),\mathbf b-\mathbf m_{\mathbf x}(s,T)\right\rangle\\
& =\left\langle\Sigma(s,t)^{-1}\mathbf y,\mathbf y\right\rangle\\
& \quad-2\left\langle\Sigma(s,t)^{-1}\mathbf y,\Sigma(s,t)\kappa(s,t)^{-1}\mathbf m_{\mathbf x}(s,t)+\Sigma(s,t)E(T,t)^\top\kappa(t,T)^{-1}\mathbf m_{\mathbf b}^-(t,T)\right\rangle\\
& \quad+\left\langle\kappa(s,t)^{-1}\mathbf m_{\mathbf x}(s,t),\mathbf m_{\mathbf x}(s,t)\right\rangle+\left\langle\kappa(t,T)^{-1}\mathbf m_{\mathbf b}^-(t,T),\mathbf m_{\mathbf b}^-(t,T)\right\rangle\\
& \quad-\left\langle\kappa(s,T)^{-1}\big(\mathbf b-\mathbf
m_{\mathbf x}(s,T)\big),\mathbf b-\mathbf m_{\mathbf
x}(s,T)\right\rangle.
\end{align*}
Since
\begin{align*}
\mathbf n_{\mathbf x,\mathbf b}(s,t) & =\Gamma(t,T)\Gamma(s,T)^{-1}\mathbf m_{\mathbf x}^+(s,t)+\Sigma(s,t)^\top\big(\Gamma(t,T)^\top\big)^{-1}\mathbf m_{\mathbf b}^-(t,T)\\
& =\Sigma(s,t)\kappa(s,t)^{-1}\mathbf m_{\mathbf x}(s,t)+\Sigma(s,t)E(T,t)^\top\kappa(t,T)^{-1}\mathbf m_{\mathbf b}^-(t,T),
\end{align*}
we have
 \begin{align*}
  \psi_{s,t}(\mathbf x,\mathbf y)
      =&\left\langle\Sigma(s,t)^{-1}\mathbf y,\mathbf y\right\rangle
        -2\left\langle\Sigma(s,t)^{-1}\mathbf y,\mathbf n_{\mathbf x,\mathbf b}(s,t)\right\rangle
        +\left\langle\kappa(s,t)^{-1}\mathbf m_{\mathbf x}(s,t),\mathbf m_{\mathbf x}(s,t)\right\rangle\\
      &+\left\langle\kappa(t,T)^{-1}\mathbf m_{\mathbf b}^-(t,T),\mathbf m_{\mathbf b}^-(t,T)\right\rangle
        -\left\langle\kappa(s,T)^{-1}\big(\mathbf b-\mathbf m_{\mathbf x}(s,T)\big),
             \mathbf b-\mathbf m_{\mathbf x}(s,T)\right\rangle.
 \end{align*}
Hence, in order to show that
$$\psi_{s,t}(\mathbf x,\mathbf y)=\left\langle\Sigma(s,t)^{-1}\big(\mathbf y-\mathbf n_{\mathbf x,\mathbf b}(s,t)\big),\mathbf y-\mathbf n_{\mathbf x,\mathbf b}(s,t)\right\rangle,$$
it remains to prove that
\begin{align*}
&\left\langle\Sigma(s,t)^{-1}\mathbf n_{\mathbf x,\mathbf b}(s,t),\mathbf n_{\mathbf x,\mathbf b}(s,t)\right\rangle\\
& \qquad=\left\langle\kappa(s,t)^{-1}\mathbf m_{\mathbf x}(s,t),\mathbf m_{\mathbf x}(s,t)\right\rangle+\left\langle\kappa(t,T)^{-1}\mathbf m_{\mathbf b}^-(t,T),\mathbf m_{\mathbf b}^-(t,T)\right\rangle\\
& \qquad\quad-\left\langle\kappa(s,T)^{-1}\big(\mathbf b-\mathbf
m_{\mathbf x}(s,T)\big),\mathbf b-\mathbf m_{\mathbf
x}(s,T)\right\rangle.
\end{align*}
The right-hand side of the above equation, due to
\begin{equation}\label{b-mx_identity}
 \mathbf b-\mathbf m_{\mathbf x}(s,T)
    =\mathbf b-E(T,s)\mathbf x-\int_s^TE(T,u)\mathbf r(u)\,{\rm d} u
    =\mathbf m_{\mathbf b}^-(t,T)-E(T,s)\mathbf m_{\mathbf x}^+(s,t),
\end{equation}
coincides with
\begin{align*}
& \left\langle \big(E(t,s)^\top\Gamma(s,t)^{-1}-E(T,s)^\top\Gamma(s,T)^{-1}\big)
               \mathbf m_{\mathbf x}^+(s,t), \mathbf m_{\mathbf x}^+(s,t)\right\rangle\\
& \qquad+\left\langle\big(\kappa(t,T)^{-1}-\kappa(s,T)^{-1}\big)\mathbf m_{\mathbf b}^-(t,T),\mathbf m_{\mathbf b}^-(t,T)\right\rangle+2\left\langle\Gamma(s,T)^{-1}\mathbf m_{\mathbf x}^+(s,t),
                   \mathbf m_{\mathbf b}^-(t,T)\right\rangle.
\end{align*}
Hence it remains to show the following three identities:
\begin{equation}\label{remain_1}\begin{split}
& E(t,s)^\top\Gamma(s,t)^{-1}-E(T,s)^\top\Gamma(s,T)^{-1}\\
& \quad=\big(\Gamma(s,T)^\top\big)^{-1}\Gamma(t,T)^\top\Sigma(s,t)^{-1}\Gamma(t,T)\Gamma(s,T)^{-1}\end{split}
\end{equation}
\begin{align}
& \kappa(t,T)^{-1}-\kappa(s,T)^{-1}
   =\Gamma(s,T)^{-1}\Gamma(s,t)\Sigma(s,t)^{-1}\Gamma(s,t)^\top\big(\Gamma(s,T)^\top\big)^{-1}\label{remain_2},\\
& \Gamma(s,T)^{-1}
   =\Gamma(s,T)^{-1}\Gamma(s,t)\Sigma(s,t)^{-1}\Gamma(t,T)\Gamma(s,T)^{-1}\label{remain_3}.
\end{align}
The validity of \eqref{remain_3} is obvious, by the definition of $\Sigma(s,t)$.
Using again the definition of $\Sigma(s,t)$, the right-hand side of \eqref{remain_2} coincides with
\begin{align*}
 \Gamma(t,T)^{-1}&\Gamma(s,t)^\top\big(\Gamma(s,T)^\top\big)^{-1}
     =\Gamma(t,T)^{-1}\Gamma(s,t)^\top\big(\Gamma(s,T)^\top\big)^{-1}\\
    &=\big(\kappa(t,T)^{-1}-\kappa(s,T)^{-1}\big)^\top
     =\kappa(t,T)^{-1}-\kappa(s,T)^{-1},
\end{align*}
where the last but one equality follows from part (a) of Lemma \ref{gen_sigmaintidentity}. Moreover, multiplying \eqref{remain_1} with $\Gamma(s,t)^\top$ from the left and $\Gamma(s,T)$ from the right, we easily obtain that \eqref{remain_1} is equivalent to
$$E(t,s)\Gamma(s,T)-\Gamma(s,t)^\top E(T,s)^\top=\Gamma(t,T),$$
and this equality holds true, by part (d) of Lemma \ref{gen_sigmaintidentity}.
\proofend

For proving almost surely continuity of the linear process bridge at the endpoint $T$,
 we recall a strong law of large numbers for continuous square integrable
 multivariate martingales with deterministic quadratic variation process due to Dzhaparidze and Spreij
 \cite[Corollary 2]{DzhSpr}; see also Koval' \cite[Corollary 1]{Kov}.
We note that the above mentioned citations are about continuous square integrable martingales with time
 interval $[0,\infty)$, but they are also valid for continuous square integrable martingales
 with time interval $[0,T)$, $T\in(0,\infty)$, with appropriate modifications in the conditions, see
 as follows (the proof of Koval' \cite[Corollary 1]{Kov} can be easily formulated for the time
 interval $[0,T)$, $T\in(0,\infty)$).

\begin{Thm}\label{THEOREM_STRONG_LAW_MULTI}
Let $T\in(0,\infty]$ be fixed and let $(\Omega,\cG,(\cG_t)_{t\in[0,T)},P)$ be a
 filtered probability space satisfying the usual conditions, i.e.,
 $(\Omega,\cG,P)$ is complete, the filtration \ $(\cG_t)_{t\in[0,T)}$ \ is right continuous,
 $\cG_0$ contains all the $P$-null sets in $\cG$ and $\cG_{T-} = \cG$, where
 $\cG_{T-}:=\sigma\left(\bigcup_{t\in[0,T)}\cG_t\right)$.
Let $(\mathbf M_t)_{t\in[0,T)}$ be an $\RR^d$-valued continuous square integrable martingale
 with respect to the filtration $(\cG_t)_{t\in[0,T)}$ such that $P(\mathbf M_0=\mathbf 0)=1$.
(The square integrability means that $\EE(m_t^{(i)})^2<\infty$, $t\in[0,T)$, $i=1,\ldots,d$,
 where $\mathbf M_t:=(m_t^{(1)},\ldots,m_t^{(d)})$, $t\in[0,T)$.)
Further, we assume that the quadratic variation process $(\langle \mathbf M\rangle_t)_{t\in[0,T)}$
 is deterministic (which yields that
 $(\langle \mathbf M\rangle_t)_{i,j}=\EE(m_t^{(i)}m_t^{(j)})$, $t\in[0,T)$, $i,j=1,\ldots,d$).
If there exists some $t_0\in[0,T)$ such that $ \langle \mathbf M\rangle_{t_0}$ is positive
 definite and $\lim_{t\uparrow T}\langle \mathbf M\rangle_t^{-1}=0\in\RR^{d\times d}$, then
 $P(\lim_{t\uparrow T}\langle \mathbf M\rangle_t^{-1}\mathbf M_t= \mathbf 0)=1$.
\end{Thm}

Using Theorem \ref{THEOREM_STRONG_LAW_MULTI} we formulate an auxiliary result which we will
 use for proving almost surely continuity of the linear process bridge at the endpoint $T$,
 see the proof of Theorem \ref{gen_multivariate_intrep_bridge}.
\begin{Lem}\label{LEMMA_almost_surely2}
Let us assume that condition \eqref{kappa_assumption} holds.
Let $T\in(0,\infty)$ be fixed and let $(\mathbf B_t)_{t\geq 0}$ be an $p$-dimensional
 standard Wiener process on a filtered probability space $(\Omega,\cA,(\cA_t)_{t\in[0,T)},P)$
 satisfying the usual conditions, constructed by the help of the standard Wiener process
 $\mathbf B$ (see, e.g., Karatzas and Shreve \cite[Section 5.2.A]{KarShr}).
The process $(\mathbf S_t)_{t\in[0,T]}$ defined by
 \begin{align}\label{SEGED_ALMOST_CONTINUITY}
    \mathbf S_t:=
         \begin{cases}
            \Gamma(t,T)
             \int_0^t \Gamma(u,T)^{-1}\Sigma(u)\,{\rm d} \mathbf B_u
               & \text{if \ $t\in[0,T)$,}\\
             \mathbf 0 & \text{if \ $t=T$,}
         \end{cases}
 \end{align}
 is a centered Gauss process with almost surely continuous paths.
\end{Lem}
\noindent{\bf Proof.}
In \eqref{SEGED_ALMOST_CONTINUITY} the stochastic integral can be understood as a usual stochastic integral
 see, e.g., Ash and Gardner \cite[Theorem 5.1.4]{AshGar}.
For this we only have to check that for all $0\leq t<T$,
$$
   \int_0^t\big((\Gamma(u,T)^{-1}\Sigma(u))_{i,j}\big)^2\,\dd u<\infty,
      \quad i=1,\ldots,d,\;\;j=1,\ldots,p,
$$
 which follows directly from part (b) of Lemma \ref{gen_sigmaintidentity}.
Due to the fact that the integrand in the stochastic integral of \eqref{SEGED_ALMOST_CONTINUITY}
 is deterministic, by Bauer \cite[Lemma 48.2]{Bau}, $(\mathbf S_t)_{t\in[0,T]}$ is a centered Gauss process.
To prove almost sure continuity, we follow the method of the proof of Lemma 5.6.9 in
 Karatzas and Sherve \cite{KarShr}.
For all $t\in[0,T)$, let
 \[
   \mathbf M_t:=\int_0^t\Gamma(u,T)^{-1}\Sigma(u)\,{\rm d} \mathbf B_u
               = \begin{pmatrix}
                   \sum_{k=1}^{p} \int_0^t (\Gamma(u,T)^{-1}\Sigma(u))_{1,k}\,{\rm d} B^{(k)}_u \\
                   \vdots \\
                   \sum_{k=1}^{p} \int_0^t (\Gamma(u,T)^{-1}\Sigma(u))_{d,k}\,{\rm d} B^{(k)}_u \\
                 \end{pmatrix},
 \]
 where $\mathbf B_t := (B_t^{(1)},\ldots,B_t^{(p)})$, $t\geq 0$.
Then, by Proposition 3.2.10 in Karatzas and Shreve \cite{KarShr},
 $(\mathbf M_t)_{t\in[0,T)}$ \ is a continuous, square integrable martingale with respect to
 the filtration $(\cA_t)_{t\in[0,T)}$ and with quadratic variation process
 $(\langle \mathbf M\rangle_t)_{t\in[0,T)}$,
 \begin{align*}
   (\langle \mathbf M\rangle_t)_{i,j}
     & = \EE\left( \sum_{k=1}^{p} \int_0^t (\Gamma(u,T)^{-1}\Sigma(u))_{i,k}\,{\rm d} B^{(k)}_u
                   \sum_{\ell=1}^{p} \int_0^t (\Gamma(u,T)^{-1}\Sigma(u))_{j,\ell}\,{\rm d}
                     B^{(\ell)}_u\right)\\
     & = \sum_{k=1}^{p}\EE\left( \int_0^t (\Gamma(u,T)^{-1}\Sigma(u))_{i,k}\,{\rm d} B^{(k)}_u
                               \int_0^t (\Gamma(u,T)^{-1}\Sigma(u))_{j,k}\,{\rm d} B^{(k)}_u\right)\\
     & = \sum_{k=1}^{p} \int_0^t (\Gamma(u,T)^{-1}\Sigma(u))_{i,k}
                              (\Gamma(u,T)^{-1}\Sigma(u))_{j,k} \,{\rm d} u \\
     & = \int_0^t \Big(\Gamma(u,T)^{-1}\Sigma(u)
                              \Sigma(u)^\top(\Gamma(u,T)^{-1})^\top\Big)_{i,j} \,{\rm d} u \\
     & = \Big(\Gamma(t,T)^{-1}\Sigma(0,t) (\Gamma(t,T)^\top)^{-1}\Big)_{i,j},
       \qquad t\in[0,T),\;\; 1\leq i,j\leq d,
 \end{align*}
 where the last equality follows by part (b) of Lemma \ref{gen_sigmaintidentity}.
Hence
 \begin{align*}
   \langle \mathbf M\rangle_t
      = \Gamma(t,T)^{-1}\Sigma(0,t) (\Gamma(t,T)^\top)^{-1},
       \qquad t\in[0,T),
 \end{align*}
 which shows that $\langle \mathbf M\rangle_t$ is symmetric and positive definite (and hence
 invertible) for all $t\in(0,T)$, by Lemma \ref{gen_sigma_identity}.

\noindent Now we check that $\lim_{t\uparrow T}\langle \mathbf M\rangle_t^{-1} = 0\in\RR^{d\times d}$.
By the definition of \ $\Sigma(0,t)$,
 \begin{align}\label{SEGED28}
   \langle \mathbf M\rangle_t
     = \Gamma(0,T)^{-1} \Gamma(0,t) (\Gamma(t,T)^\top)^{-1},\qquad t\in[0,T),
 \end{align}
 which yields that
 \begin{align*}
    & \langle \mathbf M\rangle_t^{-1}
      = \Gamma(t,T)^\top \Gamma(0,t)^{-1} \Gamma(0,T) \\
    & \;= \left(\int_t^T E(t,u)\Sigma(u)\Sigma(u)^\top E(T,u)^\top \,{\rm d} u\right)^\top
        \left(\int_0^t E(0,u)\Sigma(u)\Sigma(u)^\top E(t,u)^\top \,{\rm d} u\right)^{-1}
        \Gamma(0,T) \\
    & \; = \left(E(t,0) \int_t^T E(0,u)\Sigma(u)\Sigma(u)^\top E(0,u)^\top
                              \,{\rm d} u\, E(T,0)^\top \right)^\top \\
    & \;\phantom{=\;\;}
        \times
        \left(\int_0^t E(0,u)\Sigma(u)\Sigma(u)^\top E(0,u)^\top \,{\rm d} u\, E(t,0)^\top \right)^{-1}
        \Gamma(0,T) = 
 \end{align*}        
 \begin{align*}        
    &\; = E(T,0) \left(\int_t^T E(0,u)\Sigma(u)\Sigma(u)^\top E(0,u)^\top \,\dd u\right)^\top
        E(t,0)^\top (E(t,0)^\top)^{-1}\\
    & \; \phantom{=\;\;}
        \times\left(\int_0^t E(0,u)\Sigma(u)\Sigma(u)^\top E(0,u)^\top \,\dd u \right)^{-1}
        \Gamma(0,T),\qquad t\in(0,T).
 \end{align*}
Hence
 \begin{align*}
  \lim_{t\uparrow T}
    \langle \mathbf M\rangle_t^{-1}
       = E(T,0) \cdot \mathbf 0 \cdot
           \left(\int_0^T E(0,u)\Sigma(u)\Sigma(u)^\top E(0,u)^\top \,\dd u \right)^{-1}
           \Gamma(0,T)
           = 0\in\RR^{d\times d},
  \end{align*}
 where the inverse matrix
 \[
    \left(\int_0^T E(0,u)\Sigma(u)\Sigma(u)^\top E(0,u)^\top \,\dd u \right)^{-1}
 \]
  exists, since
 \begin{align*}
  \int_0^T E(0,u)\Sigma(u)\Sigma(u)^\top E(0,u)^\top \,\dd u
   & = E(0,T)
      \int_0^T E(T,u)\Sigma(u)\Sigma(u)^\top E(T,u)^\top \,\dd u\;
      E(0,T)^\top \\
   & = E(0,T) \kappa(0,T) E(0,T)^\top,
 \end{align*}
 and $E(0,T)$ and $\kappa(0,T)$ are invertible matrices (using also assumption \eqref{kappa_assumption}).

By Theorem \ref{THEOREM_STRONG_LAW_MULTI}, we get
 $P(\lim_{t\uparrow T} \langle \mathbf M\rangle_t^{-1} \mathbf M_t=\mathbf 0)=1$.
Then
 \begin{align*}
  \mathbf S_t
     = \Gamma(t,T) \mathbf M_t
     = \Gamma(t,T) \langle \mathbf M\rangle_t \langle \mathbf M\rangle_t^{-1}
       \mathbf M_t,\qquad t\in(0,T).
 \end{align*}
By part (a) of Lemma \ref{gen_sigmaintidentity},
 \[
    \Gamma(0,T)^{-1}\Gamma(0,t)(\Gamma(t,T)^\top)^{-1}
      = \kappa(t,T)^{-1} - \kappa(0,T)^{-1},
       \qquad 0\leq t<T,
 \]
 and hence, by \eqref{SEGED28},
 \begin{align*}
   &\Gamma(t,T)\langle \mathbf M\rangle_t
        = \Gamma(t,T) \Gamma(0,T)^{-1}\Gamma(0,t)(\Gamma(t,T)^\top)^{-1}
        = \Gamma(t,T)(\kappa(t,T)^{-1}-\kappa(0,T)^{-1}) \\
   &\; = E(t,T)\kappa(t,T)(\kappa(t,T)^{-1}-\kappa(0,T)^{-1})
       = E(t,T) - E(t,T)\kappa(t,T)\kappa(0,T)^{-1} \\
   &\; = E(t,T) - E(t,T)\int_t^T E(T,u)\Sigma(u)\Sigma(u)^\top E(T,u)^\top\,\dd u \, \kappa(0,T)^{-1} \\
   &\; = E(t,T) - E(t,T)E(T,0)\int_t^T E(0,u)\Sigma(u)\Sigma(u)^\top E(0,u)^\top\,\dd u \,
        E(T,0)^\top \kappa(0,T)^{-1},
 \end{align*}
 which yields that \ $\lim_{t\uparrow T} \Gamma(t,T)\langle \mathbf M\rangle_t=I_d$.
\ Hence we get \ $P(\lim_{t\uparrow T} \mathbf S_t=\mathbf 0)=1$.
\proofend

\vspace{1mm}

\noindent{\bf Proof of Theorem \ref{gen_multivariate_intrep_bridge}.}
Due to the fact that the integrand in the stochastic integral of \eqref{gen_multivariate_intrep}
 is deterministic, by Lemma 48.2 in Bauer \cite{Bau}, $(U_t)_{t\in[0,T)}$ \ is a Gauss process and
 the distribution of $\mathbf U_{t}$ is Gauss with mean $\EE\mathbf U_t=\mathbf n_{\mathbf a,\mathbf b}(0,t)$
 and by part (b) of Lemma \ref{gen_sigmaintidentity} we get
 $$
  \Cov(\mathbf U_{t},\mathbf U_{t})
    =\Gamma(t,T)\displaystyle\int_{0}^t\Gamma(u,T)^{-1}\Sigma(u)\Sigma(u)^\top
        \big(\Gamma(u,T)^\top\big)^{-1}\,{\rm d} u\,\Gamma(t,T)^\top=\Sigma(0,t)
 $$
for all $t\in[0,T)$. Since the deterministic functions appearing in \eqref{gen_multivariate_intrep} are
 continuous on $[0,T)$, Theorem 5.1.5 in Ash and Gardner \cite{AshGar} implies that
 $(\mathbf U_{t})_{t\in[0,T)}$ is stochastically and $L^2$-continuous.
Moreover, since the covariance-function $[0,T)\ni t\mapsto \Sigma(0,t)$ is continuous with
 $\Sigma(0,0)=0\in\RR^{d\times d}$ and $\Sigma(0,t)\to 0$ as $t\uparrow T$,
 the continuity theorem yields that $\mathbf U_t\to\mathbf b=\mathbf n_{\mathbf a,\mathbf b}(0,T)$
 in distribution as $t\uparrow T$.
Since the limit $\mathbf b$ is a constant, we get stochastic continuity
  of the extension $(\mathbf U_t)_{t\in[0,T]}$ with
 $\mathbf U_0=\mathbf a=\mathbf n_{\mathbf a,\mathbf b}(0,0)$ and
 $\mathbf U_T=\mathbf b=\mathbf n_{\mathbf a,\mathbf b}(0,T)$.
To prove the $L^2$-continuity it remains to check that $\mathbf U_t\to \mathbf b$ in $L^2$
 as $t\uparrow T$.
By part (b) of Lemma \ref{gen_sigmaintidentity}, for all $t\in[0,T)$ we get
 \begin{align*}
   \EE\Vert \mathbf U_t - \mathbf b\Vert^2
      &\leq \EE\left[\left(\Vert \mathbf n_{\mathbf a,\mathbf b}(0,t)
                     - \mathbf n_{\mathbf a,\mathbf b}(0,T)\Vert
                     + \Vert \mathbf U_t - \mathbf n_{\mathbf a,\mathbf b}(0,t)\Vert
               \right)^2\right]\\
   &\leq 2 \Vert \mathbf n_{\mathbf a,\mathbf b}(0,t) - \mathbf n_{\mathbf a,\mathbf b}(0,T)\Vert^2
        + 2\EE\left\Vert \Gamma(t,T)\int_{0}^t\Gamma(u,T)^{-1}\Sigma(u)\,{\rm d}\mathbf B_{u}\right\Vert^2 \\
   & = 2 \Vert \mathbf n_{\mathbf a,\mathbf b}(0,t) - \mathbf n_{\mathbf a,\mathbf b}(0,T)\Vert^2
       + 2\left(\sum_{i=1}^d \sum_{k=1}^{p}
         \int_{0}^t \Big(\Gamma(t,T)\Gamma(u,T)^{-1}\Sigma(u)\Big)_{i,k}^2
               \,\dd u \right) \\
   & = 2 \Vert \mathbf n_{\mathbf a,\mathbf b}(0,t) - \mathbf n_{\mathbf a,\mathbf b}(0,T)\Vert^2\\
   &\phantom{=\;}
      + 2\tr\left(\int_{0}^t \Gamma(t,T)\Gamma(u,T)^{-1}\Sigma(u)\Sigma(u)^\top (\Gamma(u,T)^{-1})^\top
                                  \Gamma(t,T)^\top\,\dd u \right) \\
   & = 2 \Vert \mathbf n_{\mathbf a,\mathbf b}(0,t) - \mathbf n_{\mathbf a,\mathbf b}(0,T)\Vert^2
      + 2\tr(\Sigma(0,t))\\
   & = 2 \Vert \mathbf n_{\mathbf a,\mathbf b}(0,t) - \mathbf n_{\mathbf a,\mathbf b}(0,T)\Vert^2
       + 2\tr\big(\Gamma(t,T)\Gamma(0,T)^{-1}\Gamma(0,t)\big)
   \to 0 \qquad \text{as \ $t\uparrow T$},
 \end{align*}
 where $\tr(A)$ denotes the trace of a squared matrix \ $A$.

Further, Lemma \ref{LEMMA_almost_surely2} yields that \ $({\mathbf U}_t)_{t\in[0,T]}$ \ is
 almost surely continuous.

Moreover, for all $0\leq s\leq t<T$ we have
\begin{align*}
\mathbf U_t = & \Gamma(t,T)\Gamma(s,T)^{-1}\\
& \times\Bigg(\phantom{\bigg|}\Gamma(s,T)\Gamma(0,T)^{-1}
      \mathbf m_{\mathbf a}^+(0,t)+\Gamma(s,T)\Gamma(t,T)^{-1}\Gamma(0,t)^\top\big(\Gamma(0,T)^\top\big)^{-1}
       \mathbf m_{\mathbf b}^-(t,T)\\
& +\Gamma(s,T)\int_{0}^t\Gamma(u,T)^{-1}\Sigma(u)\,{\rm d}\mathbf B_{u}\Bigg)\\
= & \Gamma(t,T)\Gamma(s,T)^{-1}
 \Bigg(\mathbf U_s+\Gamma(s,T)\Gamma(0,T)^{-1}\int_{s}^t E(0,u)\mathbf r(u)\,{\rm d} u\\
& \qquad \qquad \quad \quad \;\;
 +\Gamma(0,s)^\top\big(\Gamma(0,T)^\top\big)^{-1}\int_{s}^t E(T,u)\mathbf r(u)\,{\rm d} u\Bigg)\\
& +\Big(\Gamma(0,t)^\top\big(\Gamma(0,T)^\top\big)^{-1}-\Gamma(t,T)\Gamma(s,T)^{-1}\Gamma(0,s)^\top\big(\Gamma(0,T)^\top\big)^{-1}\Big)\mathbf m_{\mathbf b}^-(t,T)\\
& +\Gamma(t,T)\int_{s}^t\Gamma(u,T)^{-1}\Sigma(u)\,{\rm d}\mathbf B_{u},
\end{align*}
 where for the derivation of the last equality we used that
\begin{align*}
   &\Gamma(0,s)^\top\big(\Gamma(0,T)^\top\big)^{-1}\mathbf m_{\mathbf b}^-(s,T)
      +\Gamma(0,s)^\top\big(\Gamma(0,T)^\top\big)^{-1}
         \int_s^t E(T,u)\mathbf r(u)\,\dd u\\
   &=\Gamma(0,s)^\top\big(\Gamma(0,T)^\top\big)^{-1}
      \left(\mathbf b-\int_s^T E(T,u)\mathbf r(u)\,\dd u\right)
     +\Gamma(0,s)^\top\big(\Gamma(0,T)^\top\big)^{-1}
         \int_s^t E(T,u)\mathbf r(u)\,\dd u\\
   &=\Gamma(0,s)^\top\big(\Gamma(0,T)^\top\big)^{-1}
       \left(\mathbf b-\int_t^T E(T,u)\mathbf r(u)\,\dd u\right)\\
   &=\Gamma(0,s)^\top\big(\Gamma(0,T)^\top\big)^{-1}
      \mathbf m_{\mathbf b}^-(t,T).
\end{align*}
Next we argue that for any $0\leq s\leq u<t<T$ we have
\begin{equation}\label{gen_bridge_evol}
\Gamma(s,T)\Gamma(0,T)^{-1}E(0,u)+\Gamma(0,s)^\top\big(\Gamma(0,T)^\top\big)^{-1}E(T,u)=E(s,u).
\end{equation}
Multiplying \eqref{gen_bridge_evol} by $E(u,T)\kappa(0,T)$ from
the right, equivalently we have to show that
\begin{equation}\label{gen_bridge_evol_2}
\Gamma(s,T)+\Gamma(0,s)^\top E(T,0)^\top=E(s,T)\kappa(0,T),
\end{equation}
which holds by part (d) of Lemma \ref{gen_sigmaintidentity}.
Furthermore, by part (a) of Lemma \ref{gen_sigmaintidentity} and the symmetry of $\kappa^{-1},$ we have
\begin{align*}
\Gamma(0,t)^\top & \big(\Gamma(0,T)^\top\big)^{-1}-\Gamma(t,T)\Gamma(s,T)^{-1}\Gamma(0,s)^\top\big(\Gamma(0,T)^\top\big)^{-1}\\
& = \Gamma(t,T)\bigg[\Big(\Gamma(0,T)^{-1}\Gamma(0,t)\big(\Gamma(t,T)^\top\big)^{-1}\Big)^\top\\
& \qquad\qquad \,-\Big(\Gamma(0,T)^{-1}\Gamma(0,s)\big(\Gamma(s,T)^\top\big)^{-1}\Big)^\top\bigg]\\
& =\Gamma(t,T)\Big[\big(\kappa(t,T)^{-1}-\kappa(0,T)^{-1}\big)^\top
            -\big(\kappa(s,T)^{-1}-\kappa(0,T)^{-1}\big)^\top\Big]\\
& = \Gamma(t,T)\big(\kappa(t,T)^{-1}-\kappa(s,T)^{-1}\big)=\Gamma(s,t)^\top\big(\Gamma(s,T)^\top\big)^{-1}.
\end{align*}
Putting all together, using \eqref{gen_bridge_evol} and the above, we get
\begin{align}
\mathbf U_t= & \Gamma(t,T)\Gamma(s,T)^{-1}\mathbf m_{\mathbf U_s}^+(s,t)
 +\Gamma(s,t)^\top\big(\Gamma(s,T)^\top\big)^{-1}\mathbf m_{\mathbf b}^-(t,T)\nonumber\\
& +\Gamma(t,T)\int_{s}^t\Gamma(u,T)^{-1}\Sigma(u)\,{\rm d}\mathbf B_{u}\label{bridge_increments_solution}\\
 = & \mathbf n_{\mathbf U_s,\mathbf b}(s,t)
 +\Gamma(t,T)\int_{s}^t\Gamma(u,T)^{-1}\Sigma(u)\,{\rm d}\mathbf B_{u},\nonumber
\end{align}
where the last integral is independent of $\mathbf U_s$ (by page 434 in Bauer \cite{Bau}).
Given $\mathbf U_s=\mathbf x,$ the distribution of $\mathbf U_t$ does not depend on
 $(\mathbf U_{r})_{r\in[0,s)}$ and hence $(\mathbf U_{t})_{t\in[0,T]}$ is a Markov process.
Moreover, for any $\mathbf x\in\RR^d$ and $0\leq s<t<T$ the conditional distribution of $\mathbf U_t$ given
 $\mathbf U_s=\mathbf x$ is Gauss with mean $\mathbf n_{\mathbf x,\mathbf b}(s,t)$ and covariance matrix
 $\Sigma(s,t)$ by Lemma \ref{gen_sigmaintidentity} (b), which coincides with the Gauss density
 given by Lemma \ref{gen_multivariate_bridge}.
\proofend

\vspace{1mm}

\noindent{\bf Proof of Theorem \ref{LEMMA_multi_DE}.}
Define the $d$-dimensional Ito-process $(\mathbf V_{t})_{t\in[0,T)}$ by
$$
 \mathbf V_{t}:=\int_{0}^t\Gamma(u,T)^{-1}\Sigma(u)\,{\rm d}\mathbf B_{u},
     \quad\text{ i.e. } \quad{\rm d}
 \mathbf V_{t}=\Gamma(t,T)^{-1}\Sigma(t)\,{\rm d}\mathbf B_{t},\qquad t\in[0,T).
$$
Further, let $F:[0,T)\times\RR^d\to\RR^d$ be given by
 $F(t,\mathbf x):=\mathbf n_{\mathbf a,\mathbf b}(0,t)+\Gamma(t,T)\mathbf x$, $t\in[0,T)$, $\mathbf x\in\RR^d$.
By \eqref{gen_multivariate_intrep}, $\mathbf U_{t}=F(t,\mathbf V_{t})$ for $t\in[0,T)$.
Recall that for $0\leq s<t<T$ we have
\begin{align*}
\Gamma(s,t)  = E(s,t)\int_{s}^tE(t,u)\Sigma(u)\Sigma(u)^\top E(t,u)^\top\,{\rm d} u
             = \int_{s}^tE(s,u)\Sigma(u)\Sigma(u)^\top E(s,u)^\top\,{\rm d} u\,E(t,s)^\top,
\end{align*}
and hence using \eqref{evoderive} we get
\begin{align}\label{SEGED21}
 \begin{split}
\partial_{1}\Gamma(t,T) = & Q(t)E(t,T)\int_{t}^TE(T,u)\Sigma(u)\Sigma(u)^\top E(T,u)^\top\,{\rm d} u\\
& -E(t,T)E(T,t)\Sigma(t)\Sigma(t)^\top E(T,t)^\top\\
= & Q(t)\Gamma(t,T)-\Sigma(t)\Sigma(t)^\top E(T,t)^\top,
 \end{split}
\end{align}
and
\begin{align}\label{SEGED22}
 \begin{split}
\partial_{2}\Gamma(0,t) = & E(0,t)\Sigma(t)\Sigma(t)^\top E(0,t)^\top E(t,0)^\top\\
& +\int_{0}^tE(0,u)\Sigma(u)\Sigma(u)^\top E(0,u)^\top\,{\rm d} u\, E(t,0)^\top Q(t)^\top\\
= & \Gamma(0,t)Q(t)^\top+E(0,t)\Sigma(t)\Sigma(t)^\top.
 \end{split}
\end{align}
Further we calculate
\begin{align*}
\partial_{1}F(t,\mathbf x)= & \big(\partial_{1}\Gamma(t,T)\big)\Gamma(0,T)^{-1}
\mathbf m_{\mathbf a}^+(0,t)+\Gamma(t,T)\Gamma(0,T)^{-1}E(0,t)\mathbf r(t)\\
& +\big(\partial_{2}\Gamma(0,t)\big)^\top\big(\Gamma(0,T)^\top\big)^{-1}
\mathbf m_{\mathbf b}^-(t,T)\\
& +\Gamma(0,t)^\top\big(\Gamma(0,T)^\top\big)^{-1}E(T,t)\mathbf r(t)
+\big(\partial_{1}\Gamma(t,T)\big)\mathbf x,
\end{align*}
and $D_{2}F(t,\mathbf x)=\Gamma(t,T),$ independent of $\mathbf x\in\RR^d$. Hence, by an application of the multivariate Ito-rule, we get
$${\rm d}\mathbf U_{t}={\rm d} F(t,\mathbf V_{t})
=\partial_{1}F(t,\mathbf V_{t})\,{\rm d} t+D_{2}F(t,\mathbf V_{t})\,{\rm d}\mathbf V_{t}=\partial_{1}F(t,\mathbf V_{t})\,{\rm d} t+\Sigma(t)\,{\rm d}\mathbf B_{t},$$
where, by \eqref{SEGED21} and \eqref{SEGED22}, the coefficient of ${\rm d}t$ equals
\begin{align*}
& \partial_{1}F(t,\mathbf  V_{t})=\big(Q(t)-\Sigma(t)\Sigma(t)^\top E(T,t)^\top\Gamma(t,T)^{-1}\big)\mathbf U_{t}\\
& \quad+\Sigma(t)\Sigma(t)^\top\Big(E(0,t)^\top\big(\Gamma(0,t)^\top\big)^{-1}
+E(T,t)^\top\Gamma(t,T)^{-1}\Big)\Gamma(0,t)^\top\big(\Gamma(0,T)^\top\big)^{-1}\mathbf m_{\mathbf b}^-(t,T)\\
& \quad+\Big(\Gamma(t,T)\Gamma(0,T)^{-1}E(0,t)
+\Gamma(0,t)^\top\big(\Gamma(0,T)^\top\big)^{-1}E(T,t)\Big)\mathbf r(t),
\end{align*}
 since
\begin{align*}
  &-\Big(Q(t)-\Sigma(t)\Sigma(t)^\top E(T,t)^\top\Gamma(t,T)^{-1}\Big)
      \mathbf n_{\mathbf a,\mathbf b}(0,t)
   +(\partial_1\Gamma(t,T))\Gamma(0,T)^{-1}\mathbf m_{\mathbf a}^+(0,t)\\
  &\quad+(\partial_2\Gamma(0,t))^\top\big(\Gamma(0,T)^\top\big)^{-1}\mathbf m_{\mathbf b}^-(t,T)\\
  &=-\partial_1\Gamma(t,T) \Gamma(t,T)^{-1}\Gamma(t,T)\Gamma(0,T)^{-1} \mathbf m_{\mathbf a}^+(0,t) \\
  &\quad - \partial_1\Gamma(t,T)\Gamma(t,T)^{-1}\Gamma(0,t)^\top \big(\Gamma(0,T)^\top\big)^{-1}
                        \mathbf m_{\mathbf b}^-(t,T)\\
  &\quad +\partial_1\Gamma(t,T)\Gamma(0,T)^{-1} \mathbf m_{\mathbf a}^+(0,t)
   +(\partial_2\Gamma(0,t))^\top \big(\Gamma(0,T)^\top\big)^{-1}\mathbf m_{\mathbf b}^-(t,T)\\
  &=\Big(-\partial_1\Gamma(t,T)\Gamma(t,T)^{-1}
         +(\partial_2\Gamma(0,t))^\top\big(\Gamma(0,t)^\top\big)^{-1}\Big)
       \Gamma(0,t)^\top\big(\Gamma(0,T)^\top\big)^{-1} \mathbf m_{\mathbf b}^-(t,T)\\
  &=\Big(\Sigma(t)\Sigma(t)^\top E(T,t)^\top\Gamma(t,T)^{-1}
         -Q(t)\Gamma(t,T)\Gamma(t,T)^{-1}
         +Q(t)\Gamma(0,t)^\top \big(\Gamma(0,t)^\top\big)^{-1}\\
  &\phantom{=\;}
         +\Sigma(t)\Sigma(t)^\top E(0,t)^\top \big(\Gamma(0,t)^\top\big)^{-1}\Big)
       \Gamma(0,t)^\top\big(\Gamma(0,T)^\top\big)^{-1} \mathbf m_{\mathbf b}^-(t,T)\\
  &=\Sigma(t)\Sigma(t)^\top\Big(
          E(0,t)^\top\big(\Gamma(0,t)^\top\big)^{-1}
          +E(T,t)^\top\Gamma(t,T)^{-1}\Big)
     \Gamma(0,t)^\top\big(\Gamma(0,T)^\top\big)^{-1} \mathbf m_{\mathbf b}^-(t,T).
\end{align*}
This expression can be simplified, since by
\eqref{gen_sigma_inverse} we have
\begin{align*}
& \Big(E(0,t)^\top\big(\Gamma(0,t)^\top\big)^{-1}+E(T,t)^\top\Gamma(t,T)^{-1}\Big)\Gamma(0,t)^\top\big(\Gamma(0,T)^\top\big)^{-1}\\
& \quad=\Big(\kappa(0,t)^{-1}+E(T,t)^\top\kappa(t,T)^{-1}E(T,t)\Big)\Gamma(0,t)^\top\big(\Gamma(0,T)^\top\big)^{-1}\\
& \quad=\Sigma(0,t)^{-1}\Gamma(0,t)^\top\big(\Gamma(0,T)^\top\big)^{-1}\\
& \quad=\big(\Sigma(0,t)^\top\big)^{-1}\Gamma(0,t)^\top\big(\Gamma(0,T)^\top\big)^{-1}=\big(\Gamma(t,T)^\top\big)^{-1},
\end{align*}
and using \eqref{gen_bridge_evol_2} we further calculate
\begin{align*}
& \Gamma(t,T)\Gamma(0,T)^{-1}E(0,t)+\Gamma(0,t)^\top\big(\Gamma(0,T)^\top\big)^{-1}E(T,t)\\
& \quad=\Big(\Gamma(t,T)+\Gamma(0,t)^\top E(T,0)^\top\Big)\kappa(0,T)^{-1}E(T,t)\\
& \quad=E(t,T)E(T,t)=I_d.
\end{align*}
Putting things together, we have the SDE \eqref{gen_mult_bridge_sde}, as desired.
Since the process $(\mathbf U_t)_{t\in[0,T)}$ is adapted to the filtration $(\mathcal F_t)_{t\in[0,T)}$,
 it is a strong solution of the SDE \eqref{gen_mult_bridge_sde}.
 By Theorem 5.2.1 in {\O}ksendal \cite{Oks} or Theorem 2.32 in Chapter III
 in Jacod and Shiryaev \cite{JacShi}, strong uniqueness holds for the SDE \eqref{gen_mult_bridge_sde}.
\proofend

The following lemma is about the covariance structure of the linear process $\mathbf Z$
 and its bridge $\mathbf U$ (given in Definition \ref{DEF_bridge_multidim}).
We use this lemma in the proofs of Theorem \ref{anticipative_bridge} and Proposition
 \ref{bridge_conditioning}.
\begin{Lem}\label{covariance_identities}
For fixed $\mathbf a,\mathbf b\in\RR^d$ and $T>0$, let $(\mathbf Z_t)_{t\geq0}$ be the $d$-dimensional
 linear process given by the SDE \eqref{gen_multivariate_system} with a Gauss initial random vector
 $\mathbf Z_0$ independent of the underlying Wiener process $(\mathbf B_t)_{t\geq0}$.
 Let us suppose that condition \eqref{kappa_assumption} holds and let
 $(\mathbf U_t)_{t\in[0,T]}$ be the linear process bridge from $\mathbf a$ to $\mathbf b$ over
 $[0,T]$ derived from $\mathbf Z$ (given by Theorem \ref{gen_multivariate_intrep_bridge} and
 Definition \ref{DEF_bridge_multidim}).
Then for $0\leq s\leq t$ the covariance matrices of $\mathbf Z$ and $\mathbf U$ are given by
\begin{enumerate}
\item[(a)] $\qquad\Cov(\mathbf Z_s,\mathbf Z_t)=\Cov(\mathbf Z_t,\mathbf Z_s)^\top=\big(E(t,0)\Gamma(0,s)\big)^\top,$
\item[(b)] $\qquad\Cov(\mathbf U_s,\mathbf U_t)=\Cov(\mathbf U_t,\mathbf U_s)^\top=\big(\Gamma(t,T)\Gamma(0,T)^{-1}\Gamma(0,s)\big)^\top.$
\end{enumerate}
\end{Lem}
\noindent{\bf Proof.}
Since the stochastic integral in \eqref{increments_solution} is independent of $\mathbf Z_s,$ by \eqref{evolution} we get
\begin{align*}
\Cov(\mathbf Z_s,\mathbf Z_t) & =\Cov(\mathbf Z_s,E(t,s)\mathbf Z_s)=\Cov(\mathbf Z_s,\mathbf Z_s)E(t,s)^\top\\
& =\kappa(0,s)^\top E(0,s)^\top E(t,0)^\top=\Gamma(0,s)^\top  E(t,0)^\top,
\end{align*}
which proves (a). Similarly, using \eqref{bridge_increments_solution} and
 Theorem \ref{gen_multivariate_intrep_bridge}, we get
\begin{align*}
\Cov(\mathbf U_s,\mathbf U_t) & =\Cov(\mathbf U_s,\Gamma(t,T)\Gamma(s,T)^{-1}\mathbf U_s)=\Cov(\mathbf U_s,\mathbf U_s)\big(\Gamma(s,T)^{-1}\big)^\top\Gamma(t,T)^\top\\
 &  =\Sigma(0,s)^\top\big(\Gamma(s,T)^{-1}\big)^\top\Gamma(t,T)^\top
    =\big(\Gamma(t,T)\Gamma(0,T)^{-1}\Gamma(0,s)\big)^\top,
\end{align*}
where the last equality follows by the definition of $\Sigma(0,s).$
\proofend

\vspace{1mm}

\noindent{\bf Proof of Theorem \ref{anticipative_bridge}.}
Since the processes $\mathbf Y$ and $\mathbf U$ are almost surely continuous Gauss processes,
 using also that the law of a (continuous) stochastic process is determined
 by its finite-dimensional distributions (see, e.g., Kallenberg \cite[Proposition 2.2]{Kallenberg}),
 it is enough to show the equality of their mean and covariance functions.
Using the definition \eqref{anticipative_def} we get
\begin{align*}
 \EE\mathbf Y_t
   = \Gamma(t,T)\Gamma(0,T)^{-1}\mathbf a+\mathbf m_{\mathbf 0}(0,t)
     +\Gamma(0,t)^\top\big(\Gamma(0,T)^\top\big)^{-1}(\mathbf b-\mathbf m_{\mathbf 0}(0,T)),
\end{align*}
 for all $t\in[0,T]$.
Since
\begin{align*}
  &\mathbf m_{\mathbf 0}(0,t)
    +\Gamma(0,t)^\top\big(\Gamma(0,T)^\top\big)^{-1}
      (\mathbf b-\mathbf m_{\mathbf 0}(0,T))\\
   &=\int_0^t E(t,u)\mathbf r(u)\,\dd u
    +\Gamma(0,t)^\top \big(\Gamma(0,T)^\top\big)^{-1}
        \left(\mathbf b-\int_0^T E(T,u)\mathbf r(u)\,\dd u\right)\\
   &= E(t,0)\int_0^t E(0,u)\mathbf r(u)\,\dd u
    -\Gamma(0,t)^\top \big(\Gamma(0,T)^\top\big)^{-1}
      E(T,0)\int_0^T E(0,u)\mathbf r(u)\,\dd u\\
   &\phantom{=\;}
     +\Gamma(0,t)^\top \big(\Gamma(0,T)^\top\big)^{-1}\mathbf b\\
   &=\left(E(t,0)-\Gamma(0,t)^\top\big(\Gamma(0,T)^\top\big)^{-1}E(T,0)\right)
      \int_0^t E(0,u)\mathbf r(u)\,\dd u\\
   &\phantom{=\;}
     +\Gamma(0,t)^\top\big(\Gamma(0,T)^\top\big)^{-1}
       \left(\mathbf b-\int_t^T E(T,u)\mathbf r(u)\,\dd u\right)\\
   &=\left(E(t,0)-\Gamma(0,t)^\top\big(\Gamma(0,T)^\top\big)^{-1}E(T,0)\right)
      \int_0^t E(0,u)\mathbf r(u)\,\dd u
      +\Gamma(0,t)^\top\big(\Gamma(0,T)^\top\big)^{-1}
       \mathbf m_{\mathbf b}^-(t,T),
\end{align*}
 we get
\begin{align*}
\EE\mathbf Y_t = & \Gamma(t,T)\Gamma(0,T)^{-1}\\
 & \times\left(\mathbf a+\Gamma(0,T)\Gamma(t,T)^{-1}
    \left[E(t,0)-\Gamma(0,t)^\top\big(\Gamma(0,T)^\top\big)^{-1}E(T,0)\right]
    \int_0^tE(0,u)\mathbf r(u)\,\dd u\right)\\
& +\Gamma(0,t)^\top\big(\Gamma(0,T)^\top\big)^{-1}\mathbf m_{\mathbf b}^-(t,T)\\
= & \mathbf n_{\mathbf a,\mathbf b}(0,t)
= \EE \mathbf U_t,\quad t\in[0,T],
\end{align*}
where the last but one equality follows by part (c) of Lemma \ref{gen_sigmaintidentity}
together with \eqref{bridge_mean}, and the last equality by Theorem \ref{gen_multivariate_intrep_bridge}.
Further, using part (a) of Lemma \ref{covariance_identities}, for all $0\leq s\leq t\leq T$ we have
 \begin{align*}
 \Cov(\mathbf Y_s,\mathbf Y_t)
    = & \Cov\Big(\mathbf Z_s-\Gamma(0,s)^\top\big(\Gamma(0,T)^\top\big)^{-1}\mathbf Z_T,
                \mathbf Z_t-\Gamma(0,t)^\top\big(\Gamma(0,T)^\top\big)^{-1}\mathbf Z_T\Big)\\
 = & \Cov(\mathbf Z_s,\mathbf Z_t)-\Gamma(0,s)^\top\big(\Gamma(0,T)^\top\big)^{-1}\Cov(\mathbf Z_T,\mathbf Z_t)\\
  & -\Cov(\mathbf Z_s,\mathbf Z_T)\Gamma(0,T)^{-1}\Gamma(0,t)\\
  & +\Gamma(0,s)^\top\big(\Gamma(0,T)^\top\big)^{-1}\Cov(\mathbf Z_T,\mathbf Z_T)\Gamma(0,T)^{-1}\Gamma(0,t)\\
 = & \Gamma(0,s)^\top E(t,0)^\top-\Gamma(0,s)^\top\big(\Gamma(0,T)^\top\big)^{-1}E(T,0)\Gamma(0,t)\\
 & -\Gamma(0,s)^\top E(T,0)^\top\Gamma(0,T)^{-1}\Gamma(0,t)\\
 & +\Gamma(0,s)^\top\big(\Gamma(0,T)^\top\big)^{-1}\kappa(0,T)\Gamma(0,T)^{-1}\Gamma(0,t),
 \end{align*}
 and hence
 \begin{align*}
  \Cov(\mathbf Y_s,\mathbf Y_t)
   = & \Gamma(0,s)^\top\big(\Gamma(0,T)^\top\big)^{-1}\left[\Gamma(0,T)^\top E(t,0)^\top-E(T,0)\Gamma(0,t)\right.\\
     & \qquad\left.-\Gamma(0,T)^\top E(T,0)^\top\Gamma(0,T)^{-1}\Gamma(0,t)+\kappa(0,T)\Gamma(0,T)^{-1}\Gamma(0,t)\right]\\
   = & \Gamma(0,s)^\top\big(\Gamma(0,T)^\top\big)^{-1}\left[\Gamma(0,T)^\top E(t,0)^\top-E(T,0)\Gamma(0,t)\right]\\
   = & \Gamma(0,s)^\top\big(\Gamma(0,T)^\top\big)^{-1}\Gamma(t,T)^\top,
 \end{align*}
where the last equality follows by part (d) of Lemma \ref{gen_sigmaintidentity}.
The last line of the above equation coincides with $\Cov(\mathbf U_s,\mathbf U_t)$
by part (b) of Lemma \ref{covariance_identities}, as desired.
\proofend

\vspace{1mm}

\noindent{\bf Proof of Proposition \ref{bridge_conditioning}.}
Since $(\mathbf Z_t)_{t\geq 0}$ is a Gauss process, by Shiryaev \cite[Theorem 2, page 303]{Shi},
 the conditional distribution of $\mathbf Z:=(\mathbf Z_{t_1}^\top,\ldots,\mathbf Z_{t_n}^\top)^\top$ given
 $\mathbf Z_T=\mathbf b$ is known to be a $\RR^{dn}$-dimensional Gauss distribution with mean vector
 $$
 \mathbf m:=\EE\mathbf Z+\left[\bigotimes_{i=1}^n\Cov(\mathbf Z_{t_i},\mathbf Z_T)\right]
             \Cov(\mathbf Z_T,\mathbf Z_T)^{-1}\big(\mathbf b-\EE\mathbf Z_T\big),
 $$
 and with covariance matrix
 $$
 C:=\Cov(\mathbf Z,\mathbf Z)-\left[\bigotimes_{i=1}^n\Cov(\mathbf Z_{t_i},\mathbf Z_T)\right]
    \Cov(\mathbf Z_T,\mathbf Z_T)^{-1}\left[\bigotimes_{i=1}^n\Cov(\mathbf Z_{t_i},\mathbf Z_T)\right]^\top,
 $$
 where, due to symmetry, the $(nd\times d)$-matrix $\bigotimes_{i=1}^n\Cov(\mathbf Z_{t_i},\mathbf Z_T)$
 is defined by
$$
 \bigotimes_{i=1}^n\Cov(\mathbf Z_{t_i},\mathbf Z_T)
   :=\big(\Cov(\mathbf Z_{t_1},\mathbf Z_T),\ldots,\Cov(\mathbf Z_{t_n},\mathbf Z_T)\big)^\top.
$$
Due to the fact that the integrand in the stochastic integral of \eqref{gen_multivariate_intrep} is
 deterministic, by Lemma 48.2 in Bauer \cite{Bau}, the process  $(\mathbf U_t)_{t\in[0,T]}$ is a Gauss process,
 and hence $(\mathbf U_{t_1}^\top,\ldots,\mathbf U_{t_n}^\top)^\top$ is a $nd$-dimensional Gauss distributed
 random variable.
Then all we have to do is to check that the mean vector and the covariance matrix of
 $(\mathbf U_{t_1}^\top,\ldots,\mathbf U_{t_n}^\top)^\top$ coincides with the mean vector and
 the covariance matrix given above.
Using Theorem \ref{gen_multivariate_intrep_bridge}, due to the tensor calculus it is sufficient to prove that
 for all $0\leq i\leq j\leq n$
\begin{equation}\label{mean_conditioning}
\mathbf m_i:=\EE\mathbf Z_{t_i}+\Cov(\mathbf Z_{t_i},\mathbf Z_T)\Cov(\mathbf Z_T,\mathbf Z_T)^{-1}\big(\mathbf b-\EE\mathbf Z_T\big)=\EE\mathbf U_{t_i}=\mathbf  n_{\mathbf a,\mathbf b}(0,t_i),
\end{equation}
and
\begin{equation}\label{covariance_conditioning}\begin{split}
C_{ij} & :=\Cov(\mathbf Z_{t_i},\mathbf Z_{t_j})-\Cov(\mathbf Z_{t_i},\mathbf Z_T)\Cov(\mathbf Z_T,\mathbf Z_T)^{-1}\Cov(\mathbf Z_{t_j},\mathbf Z_T)^\top\\
& = \Cov(\mathbf U_{t_i},\mathbf U_{t_j})=\big(\Gamma(t_j,T)\Gamma(0,T)^{-1}\Gamma(0,t_i)\big)^\top,
\end{split}\end{equation}
where the last equality follows from part (b) of Lemma \ref{covariance_identities}. Using part (a) of Lemma \ref{covariance_identities} and \eqref{b-mx_identity} we calculate
\begin{align*}
\mathbf m_i & =\mathbf m_{\mathbf a}(0,t_i)+\Gamma(0,t_i)^\top E(T,0)^\top\big(\Gamma(0,T)^\top E(T,0)^\top\big)^{-1}\big(\mathbf b-\mathbf m_{\mathbf a}(0,T)\big)\\
& =E(t_i,0)\mathbf m_{\mathbf a}^+(0,t_i)+\Gamma(0,t_i)^\top \big(\Gamma(0,T)^\top \big)^{-1}\big(\mathbf m_{\mathbf b}^-(t_i,T)-E(T,0)\mathbf m_{\mathbf a}^+(0,t_i)\big)\\
& =\Gamma(t_i,T)\Gamma(0,T)^{-1}\mathbf m_{\mathbf a}^+(0,t_i)+\Gamma(0,t_i)^\top \big(\Gamma(0,T)^\top \big)^{-1}\mathbf m_{\mathbf b}^-(t_i,T),
\end{align*}
 where the last equality follows by part (c) of Lemma \ref{gen_sigmaintidentity}. This shows \eqref{mean_conditioning} by using \eqref{bridge_mean}. Moreover, using part (a) of Lemma \ref{covariance_identities} we have
\begin{align*}
C_{ij} & =\Gamma(0,t_i)^\top E(t_j,0)^\top-\Gamma(0,t_i)^\top E(T,0)^\top\big(\Gamma(0,T)^\top E(T,0)^\top\big)^{-1}E(T,0)\Gamma(0,t_j)\\
& =\Gamma(0,t_i)^\top\big(\Gamma(0,T)^\top\big)^{-1}\left[\Gamma(0,T)^\top E(t_j,0)^\top-E(T,0)\Gamma(0,t_j)\right]\\
& =\big(\Gamma(t_j,T)\Gamma(0,T)^{-1}\Gamma(0,t_i)\big)^\top,
\end{align*}
where the last equality follows by part (d) of Lemma \ref{gen_sigmaintidentity}.
This shows \eqref{covariance_conditioning}.
\proofend

\vspace{1mm}

\noindent{\bf Proof of Theorem \ref{THEOREM2}.}
Using the notations of Section \ref{SECTION_MULTI}, we have
 \begin{align*}
   &\Phi(t) = E(t,0) = \ee^{\overline q(t)},\qquad t\geq 0,\\
   &E(t,s) = E(t,0)E(0,s)= E(t,0)E(s,0)^{-1}
           =\ee^{\overline q(t)-\overline q(s)},
           \qquad t,s\geq 0,
 \end{align*}
 and hence for all \ $0\leq s<t$, \ we get \ $\kappa(s,t)=\gamma(s,t)$ \ and
 \begin{align*}
    &\Gamma(s,t) = \ee^{\overline q(s)-\overline q(t)}\gamma(s,t),\\[1mm]
    &\Sigma(s,t)
       =\ee^{\overline q(t)-\overline q(T)}\gamma(t,T)
        \big(\ee^{\overline q(s)-\overline q(T)}\gamma(s,T)\big)^{-1}
        \ee^{\overline q(s)-\overline q(t)}\gamma(s,t)
       =\frac{\gamma(s,t)\gamma(t,T)}{\gamma(s,T)}
       =\sigma(s,t).
 \end{align*}
Then for all \ $0\leq s <t<T$,
 \begin{align*}
   &\Gamma(t,T)\Gamma(s,T)^{-1}m_a^+(s,t) + \Gamma(s,t)\Gamma(s,T)^{-1}m_b^-(t,T)  \\
     & = \frac{\ee^{\overline q(t)-\overline q(T)}\gamma(t,T)}
             {\ee^{\overline q(s)-\overline q(T)}\gamma(s,T)}
        \left(a+\int_s^t E(s,u)r(u)\,\dd u\right)
        +\frac{\ee^{\overline q(s)-\overline q(t)}\gamma(s,t)}
             {\ee^{\overline q(s)-\overline q(T)}\gamma(s,T)}
         \left(b-\int_t^T E(T,u)r(u)\,\dd u\right) \\
     & = \ee^{\overline q(t)-\overline q(s)}
        \frac{\gamma(t,T)}{\gamma(s,T)}
        \left(a+\int_s^t \ee^{\overline q(s)-\overline q(u)} r(u)\,\dd u\right)
        + \ee^{\overline q(T)-\overline q(t)}
         \frac{\gamma(s,t)}{\gamma(s,T)}
         \left(b-\int_t^T \ee^{\overline q(T)-\overline q(u)} r(u)\,\dd u\right) \\
      & = \frac{\gamma(t,T)}{\gamma(s,T)} m_a(s,t)
          + \ee^{\overline q(T)-\overline q(t)}
         \frac{\gamma(s,t)}{\gamma(s,T)}
         \left(b-\int_t^T \ee^{\overline q(T)-\overline q(u)} r(u)\,\dd u\right)
       =n_a(s,t).
 \end{align*}
Hence Theorem \ref{gen_multivariate_intrep_bridge} yields that for all $t\in[0,T)$,
 \begin{align*}
   U_t &= n_{a,b}(0,t)
         + \ee^{\overline q(t) - \overline q(T)} \gamma(t,T)
           \int_0^t \big(\ee^{\overline q(u) - \overline q(T)} \gamma(u,T)\big)^{-1}\sigma(u)\,\dd B_u \\
       &= n_{a,b}(0,t)
         + \int_0^t \frac{\gamma(t,T)}{\gamma(u,T)}
                    \ee^{\overline q(t) - \overline q(u)}\sigma(u)\,\dd B_u ,
 \end{align*}
 as desired.
\proofend

\vspace{1mm}

\noindent{\bf Proof of Proposition \ref{LEMMA5}.}
 It is enough to check that the process \ $(Z_t^*)_{t\geq 0}$ \
 has the same finite dimensional distributions as
 $$
   Z_t=m_0(0,t)+\int_0^t\ee^{\bar q(t)-\bar q(s)}\sigma(s)\,\dd B_s,
        \quad t\geq 0,
 $$
has. Indeed, \ $(Z_t)_{t\geq 0}$ \ is a strong solution of the SDE
 \eqref{gen_OU_egyenlet} and the law of a stochastic process is determined by its finite
 dimensional distributions (see, e.g., Kallenberg \cite[Proposition 2.2]{Kallenberg}).
Hence it is sufficient to check that \ $(Z_t^*)_{t\geq 0}$ \ is a Gauss process
 with the same expectation and covariance function that \ $Z$ \ has. Since \ $B^*$ \ is a Gauss
process we have \ $Z^*$ \ is also a Gauss process. The expectation
functions are the same, since
 \ $\EE Z_t=\EE Z_t^*=m_0(0,t),$ $t\geq 0.$ \ Moreover, using \eqref{SEGED26} one can calculate
 $$
   \Cov(Z_s,Z_t)
    =\ee^{\bar q(s)+\bar q(t)}\int_0^{s\wedge t}\ee^{-2\bar q(u)}\sigma^2(u)\,\dd u,
      \quad s,t\geq 0.
 $$
  One can also easily derive that
 \begin{align*}
   \Cov(Z_s^*,Z_t^*)
      &=\ee^{\bar q(s)+\bar q(t)}\Cov\left(B^*(\ee^{-2\bar q(s)}\gamma(0,s)),B^*(\ee^{-2\bar q(t)}\gamma(0,t))\right)\\
      &=\ee^{\bar q(s)+\bar q(t)}\min(\ee^{-2\bar q(s)}\gamma(0,s),\ee^{-2\bar q(t)}\gamma(0,t)),
         \quad s,t\geq 0.
 \end{align*}
 Since the function
 $$
   [0,\infty)\ni t\mapsto \ee^{-2\bar q(t)}\gamma(0,t)
       =\int_0^t\ee^{-2\bar q(u)}\sigma^2(u)\;\dd u,
 $$
  is monotone increasing, we have
  $$
   \Cov(Z_s^*,Z_t^*)
    =\ee^{\bar q(s)+\bar q(t)}\int_0^{s\wedge t}\ee^{-2\bar q(u)}\sigma^2(u)\,\dd u,
      \quad s,t\geq 0,
  $$
as desired.
\proofend

\end{document}